\documentclass[11pt]{article}

\usepackage[utf8]{inputenc}
\usepackage[a4paper, total={7in, 9.5in}]{geometry}
\usepackage{amsmath,amssymb,amsthm,bbm,graphicx,cases,float,wrapfig,multicol,epsfig,float}
\usepackage{amsfonts,amssymb,mathrsfs,color,url,amsthm,bm,mathtools}

\usepackage{tikz}
\usepackage{hyperref}
\theoremstyle{thmstyleone}%
\newtheorem{theorem}{Theorem}%  meant for continuous numbers
\newtheorem{corollary}{Corollary}%
\newtheorem{proposition}{Proposition}% to get separate numbers for theorem and proposition etc.
\theoremstyle{thmstyletwo}%
\newtheorem{example}{Example}%
\newtheorem{remark}{Remark}%
\theoremstyle{thmstylethree}%

\usepackage[doublespacing]{setspace}
\usepackage[fontsize=12pt]{fontsize}
\usepackage{booktabs}

\newcommand{\sftwo}[4]{ #1 \left[\scriptsize{\begin{matrix*}[c]  
                                                        #2 \\ 
                                                        #3
                                                  \end{matrix*}}  \; #4 \right]}
                                                  
\newcommand{\sffour}[6]{ #1 \small{\left[\scriptsize{\begin{matrix*}[c]  
                                                        #2  &  #3    \\ 
                                                        #4  &  #5 
                                                  \end{matrix*}}  \; #6 \right]}}
\newcommand{\sfsix}[8]{ #1 \small{\left[\scriptsize{\begin{matrix*}[c]  
                                                        #2  &  #3  &  #4  \\ 
                                                        #5  &  #6  &  #7
                                                  \end{matrix*}}  \; #8 \right]}}

\newcommand{\bI}{\mathbf{I}}

\newcommand{\bA}{\mathbf{A}}

\newcommand{\bX}{\mathbf{X}}

\newcommand{\bbeta}{\boldsymbol\theta}

\newcommand{\bS}{\mathbf{S}}

\newcommand{\by}{\mathbf{y}}

\newcommand{\R}{\mathbb R}

\newcommand{\N}{\mathbb N}

\newcommand{\id}{\,\mathrm{d}}

\newcommand{\dd}{\text{d}}
\DeclareMathOperator*{\argmin}{argmin} 
\DeclareMathOperator*{\argmax}{argmax}

\DeclareMathOperator{\diag}{diag}
\DeclareMathOperator{\pen}{pen}

\DeclareMathOperator{\Beta}{B}

\DeclareMathOperator{\GG}{GG}
\DeclareMathOperator{\TPB}{G3B}
\DeclareMathOperator{\EP}{EP}

\DeclareMathOperator{\Gambel}{Gambel}

\DeclareMathOperator{\Prob}{\mathbb{P}}

\usepackage{natbib}
\usepackage[noblocks]{authblk}

\title{Singularities in Bayesian Inference: Crucial or Overstated?}
\author[1]{Maria De Iorio} 
\author[2]{Andreas Heinecke}
\author[3]{Beatrice Franzolini}
\author[4]{Rafael Cabral}
\affil[1,4]{Yong Loo Lin School of Medicine, National University of Singapore,  Singapore}
\affil[2]{Yale-NUS College,  Singapore}
\affil[1]{Institute for Human Development and Potential, A*STAR}
\affil[3]{Bocconi Institute for Data Science and Analytics, Bocconi University, Italy}
\date{}

\begin{document}

\maketitle

\abstract{Over the past two decades, shrinkage priors have become increasingly popular, and many proposals can be found in the literature. 
These priors aim to shrink small effects to zero while maintaining true large effects. Horseshoe-type priors have been particularly successful in various applications, mainly due to their computational advantages. 
However, there is no clear guidance on choosing the most appropriate prior for a specific setting. 
In this work, we propose a framework that encompasses a large class of shrinkage distributions, including priors with and without a singularity at zero. 
By reframing such priors in the context of reliability theory and wealth distributions, we provide insights into the prior parameters and shrinkage properties. 
The paper's key contributions are based on studying the folded version of such distributions, which we refer to as the Gambel distribution. 
The Gambel can be rewritten as the ratio between a Generalised Gamma and a Generalised Beta of the second kind. 
This representation allows us to gain insights into the behaviours near the origin and along the tails, compute measures to compare their distributional properties, derive consistency results, devise MCMC schemes for posterior inference and ultimately provide guidance on the choice of the hyperparameters. }

\section{Introduction}
In a variety of fields, ranging from signal processing to statistics and machine learning, extensive attention has been devoted to the problem of estimating a sparse vector of regression coefficients.
As a main example, consider the classical linear regression problem of estimating an unknown signal $\bbeta = (\theta_1, \ldots, \theta_p)\in\R^p$ from an $n$-dimensional observation vector $\mathbf{y}\in\R^n$, assumed to be generated according to  
$\mathbf{y}= \mathbf{X} \boldsymbol\theta +\boldsymbol\epsilon$
where $\bX\in\R^{n\times p}$ is known and $\boldsymbol\epsilon\sim \mathcal{N}_n(\mathbf{0},\sigma^2\bI_n)$ follows a multivariate normal distribution with independent and homoscedastic components.
In many high-dimensional problems (large $p$), it can be assumed that the signal generating the observation is either sparse or compressible. 
Precisely, the \emph{true} signal $\boldsymbol\theta$ is called \emph{$s$-sparse} if at most $s\leq p$ of its entries are non-zero, and \emph{compressible} when the ordered absolute values of its entries rapidly decrease in amplitudes and, thus, can be well approximated by a sparse vector \citep[see, e.g.,][]{dziwoki2020sparse}. 
The problem of estimating sparse or compressible signals is routinely approached via regularised least squares estimators.
Assuming without loss of generality $\sigma^2=1$, a regularised least squares estimator for $\boldsymbol\theta$ is obtained as
\begin{equation}
    \label{Reg_LS}
\widehat{\boldsymbol\theta} =
\argmin_{\boldsymbol\theta\in\R^p} \left\{\tfrac{1}{2} \| \by-\bX\boldsymbol\theta\|^2 + \pen(\boldsymbol\theta)\right\}
\end{equation}
under a sparsity-promoting and separable (i.e., additive) penalty function $\pen(\boldsymbol\theta)\colon \R^p\to \R$. 
Famous examples of penalty functions are Lasso \citep{tibshirani1996regression},
i.e., $\pen(\boldsymbol\theta) = \lambda\sum_{j=1}^p|\theta_j|$; 
and Elastic-Net \citep{zou2005regularization},
i.e., $\pen(\boldsymbol\theta)=\lambda_1\sum_{j=1}^p|\theta_j|+\lambda_2\sum_{j=1}^p\theta_j^2$. 
Regularised estimators come with many useful advantages over non-regularised estimators (like ordinary least squares), such as the possibility of performing variable selection, reframing of the inferential problem into a well-defined optimisation problem also for $p>n$ scenarios, and exploiting of the bias-variance trade-off to improve out-of-sample prediction performance of regression models. For detailed accounts on regularised estimators for variable selection, see, for instance, \cite{lv2009unified} and \cite{wen2018survey}. 
The connection between penalised/regularised regression and Bayesian shrinkage priors is well known. 
To favour a sparse solution in the Bayesian framework, a prior $\pi(\bbeta)$ placed on the regression coefficients aims at shrinking small entries in the signal towards zero while maintaining ``true'' large effects. 
The posterior distribution corresponding to the regression model and a prior $\pi(\bbeta)$ for the coefficients is  $\pi(\bbeta\mid\by)\propto  \pi(\by\mid\bbeta)\pi(\bbeta)$, where $\pi(\by\mid\bbeta)$ is the density of a $n$-variate Normal distribution with expectation $\bX\boldsymbol\theta$ and covariance matrix $\sigma^2\bI_n$, resulting from the likelihood. Assuming $\sigma^2=1$, the Bayesian maximum-a-posteriori (MAP) estimator is then %for the discrete $0$-$1$-loss function
\begin{align*}
\argmax_{\bbeta\in\R^p}\pi(\bbeta\mid\by)
=\argmin_{\bbeta\in\R^p} \left\{\tfrac{1}{2}\|\by- \bX \bbeta \|^2-\log \pi(\bbeta)\right\}
\end{align*} 
which is equivalent to~\eqref{Reg_LS} for 
 $\pen(\bbeta)= - \log \pi(\bbeta)$. 
 Moreover, $\pi(\bbeta)$ is often defined via (conditional) independent priors for the components of the signal, i.e., $\theta_j \mid \eta \overset{iid}{\sim} \pi (\theta_j; \eta)$, for $j=1, \ldots, p$, for some auxiliary parameter $\eta$. Conditional independence (given $\eta$) in the Bayesian framework mimics the separability of the corresponding penalty in \eqref{Reg_LS}. 
For instance, the popular Lasso is equivalent to the MAP estimator obtained as  $\pi(\bbeta) = \prod_{j=1}^p (2 \lambda)^{-1}\exp\{-\lambda|\theta_j|\}$ from independent and identical Laplace priors for each coefficient.
 
 A wealth of shrinkage priors have been proposed in the Bayesian literature \citep[see, e.g.,][]{hsiang1975bayesian,carvalho2009handling,polson2009alternative,brown2010inference,HorseshoeB,armagan2011generalized,armagan2013posterior,polson2014bayesian,bhattacharya2015dirichlet,zhang2016high,bhadra2017horseshoe2,bai2018beta,rovckova2018spike}, often with little theoretical justification for the comparative performance of these priors and limited guidance on hyperparameter selection.    
 The choice of a shrinkage prior in a particular sparse regression application is usually driven by qualitative and computational reasoning. 
Arguably, this is due to a lack of criteria or measures to compare the shrinkage behaviour of such priors. 
Such problems are further amplified 
by the fact that shrinkage priors are usually obtained as scale mixture distributions, leading to a marginal distribution of the regression coefficients corresponding to a higher transcendental function. 
Considerations to distinguish among different amounts of shrinkage are often posed on the behaviour of the prior around zero \citep[see, for instance,][]{carvalho2009handling, bhadra2019lasso} and, in particular, we may distinguish between priors with a point mass at zero \citep[i.e., two-groups models or Spike-and-Slab priors associating positive probability to the event $\theta_j = 0 $; see, e.g.,][]{mitchell1988bayesian,george2000variable,george2000calibration,yuan2005efficient} and continuous priors which associate null probability to the event $\theta_j = 0 $ \citep[see, for instance][]{armagan2013generalized,carvalho2010horseshoe,griffin2010inference, hans2011elastic,park2008bayesian,zhang2012ep}.   
Among continuous priors, we can further distinguish between priors with unbounded density at zero (to which we refer as shrinkage priors \emph{with singularity}) or continuous bounded density (to which we refer as shrinkage priors  \emph{without singurality}). Nonetheless, a comprehensive and general study of the distributional properties of shrinkage prior distributions is currently lacking, though it can provide useful insights into their shrinkage behaviour. 
For instance, we show how continuous priors with and without a singularity at the origin often perform similarly in terms of shrinkage if other characteristics, e.g., tail shape and Gini index, are similar. 
In this regard, also \cite{song2023nearly} - considering general shrinkage priors - obtained asymptotic rates of convergence that depend not only on the mass of the prior around the origin but also on the fatness or flatness of the tails. 
In particular, considering polynomial-tailed priors, they conclude that under certain conditions the posterior contraction rate and the variable selection performance of continuous shrinkage priors are close to those of priors with a point mass.

To build a general framework encompassing a variety of shrinkage priors and to study their distributional properties, we draw a connection between the problem of sparse recovery using shrinkage priors (or sparsity-inducing penalties) studied in statistics as well as machine learning, and the theory of rare/extreme events, wealth distributions, and reliability, widely studied in the finance/actuarial and engineering literature. 
Such a connection should not be surprising, noting that there exists a substantial similarity between absolute-valued entries of sparse or compressible vectors and random samples generated from \emph{size distributions}, i.e., right-asymmetric and fat-tailed distributions originally used to model income and wealth \citep{kleiber2003statistical}. 
However, in contrast to the reliability and economic theory literature, the literature about shrinkage priors still requires a proper quantification of the behaviour of the tails, of the hazard function, and the overall inequality indexes of most used shrinkage distributions. Indeed, the tail behaviour of shrinkage priors has generally been overlooked, and researchers have often focussed on behaviour at the origin. 
Our main contributions are obtained by exploiting this connection, leveraging results and concepts from wealth distributions and reliability, and applying them to sparse recovery.
This approach suggests many interesting quantitative directions for comparing different shrinkage priors and for studying the sparsity level induced by different shrinkage priors, such as the assessment of hazard behaviour, tail indexes, Gini indexes, and Lorenz ordering, which are all considered in the following sections and which can be widely applied to various distributions.

\section{Preliminaries}\label{sec:preliminaries}
A popular way to devise a shrinkage prior is by fattening the tails of a Normal distribution via stochastic volatility \citep{bhattacharya2012bayesian}, i.e., via a Scale-Mixture-of-Normals (SMN) obtained by defining the prior $\pi(\theta)$ for each entry $\theta$ of the signal $\bbeta$ hierarchically as 
$\theta\mid\lambda \sim \mathcal{N}(0, \lambda^2)$ and $\lambda\sim g(\lambda)$,
where $\mathcal{N}(0, \lambda^2)$ denotes a Normal distribution with zero mean and variance $\lambda^2$, and where $g$ is a continuous distribution supported on the positive real line. 
When this is the case, as in the Horseshoe (HS) prior literature \citep{carvalho2009handling, HorseshoeB, van2014horseshoe,ghosh2015asymptotic,van2017adaptive,piironen2017sparsity,bhadra2017horseshoe,bhadra2019lasso}, the parameter $\kappa=1/(1+\lambda^2)$ is referred to as \emph{shrinkage weight}, while $\lambda^2$ is called \emph{shrinkage coefficient}. 
The shrinkage weight is interpreted as a random shrinkage parameter, while the \emph{significance weight}
$1-\kappa$ is considered as an inclusion probability, which can be interpreted similarly to the probability mass at zero of a Spike-and-Slab (SnS) prior. More generally, we say that a shrinkage prior  $\pi$ is a scale mixture whenever 
\begin{align} 
\pi(\theta) = \int_0^\infty f(\theta| \lambda) g(\lambda) \id \lambda
\label{eq:mixprior}
\end{align}
where $g$ is a continuous distribution and $\lambda$ is a scale parameter of the density $f$.
Many proposals have been made for the choice of a continuous prior on the variance in \eqref{eq:mixprior} \citep[see, e.g.,][]{tiao1965bayesian,gelman2006,morris2011estimating}. Roughly speaking, the literature on scale mixtures can be divided into two classes: mixing with (generalised) Gamma/inverse-Gamma or mixing with (generalised) Beta/inverse-Beta. \cite{zhang2012ep} provide a unifying framework for the first case, introducing a family of sparsity-inducing priors called Exponential Power-Generalized Inverse Gaussian (EP-GIG). 
In this work, we focus on the second case (to which the HS-like priors belong), while most of the criteria we devise to study their distributional properties could also be extended to other classes of priors.
The class of shrinkage priors presented and studied in the following sections can be represented by mixing the shrinkage weight of an Exponential Power (EP) distribution with a generalised three-parameter Beta (G3B) density. 
Even though the strategy to create the mixture is different, our class of priors and EP-GIG priors overlap, as detailed later. 
The main results obtained in this work arrive from the study of the folded version of the EP-G3B mixtures on $\mathbb{R}^+$, which we refer to as the \emph{Gambel distribution}, since it can be expressed as the ratio between a Generalised Gamma (GG) distribution and a Generalised Beta distribution of the second kind (GB2). 
Thus, four main classes of distributions are strictly connected to Gambel priors: Exponential Power (EP), GB2, Generalized Gamma (GG), and G3B. 
The probability density functions and parameter ranges for these distributions are summarized in Table~\ref{tab:pdfs}, where $\Gamma(x)$ denotes the Gamma function, $\phi(q) = \sqrt{\Gamma(3/q)/\Gamma(1/q)}$, and $B(a,b)$ denotes the Beta function.
\begin{table}[tbh]
\centering
\caption{Probability density functions for various distributions.}
\resizebox{\textwidth}{!}{
\begin{tabular}{l|l|l}
\toprule
 Probability Density Function & Support & Parameters\\
\midrule
$ \displaystyle \text{EP}(x \mid q,\kappa) = \frac{q\,\phi(q)\exp\left((-\kappa (\phi(q)\, |x|)^q)/(1-\kappa)\right)}{2(\frac{1}{\kappa} - 1)^{1/q} \Gamma(1/q)}$ & $\mathbb{R}$ & $q >0, \kappa \in (0,1)$ \\%[0.25cm]
\midrule
$ \displaystyle \text{GB2}(x\mid p, q, a, b) = \frac{p (x/q)^{ap-1} (1+(x/q)^p)^{-(a+b)}}{q B(a,b)}$ & $\mathbb{R}^+$ & $p>0$, $q>0$, $a>0$, $b>0$ \\%[0.25cm]
\midrule
$ \displaystyle \text{GG}(x \mid a, d, q) = \frac{q}{a^d \Gamma(d/q)}x^{d-1}\exp\left(-\left(\frac{x}{a}\right)^q\right)$ &  $\mathbb{R}^+$ & $a>0$, $d>0$, $q>0$ \\%[0.25cm]
\midrule
$ \displaystyle \text{G3B}(x \mid a, b, \xi) = \frac{1}{B(a,b)} \frac{\xi^a x^{a-1} (1-x)^{b-1}}{(1-(1-\xi)x)^{a+b}}$ & $(0,1)$ & $a>0$, $b>0$, $\xi>0$ \\%[0.25cm]
\bottomrule
\end{tabular}}
\label{tab:pdfs}
\end{table}

EP distributions have zero mean, and standard deviation 
$\lambda : =\left(\frac{1-\kappa}{\kappa}\right)^{1/q}$.
To the best of our knowledge, the  EP has first been described in \cite{subbotin1923law} and has then been studied in numerous works \citep[see, e.g., ][for a comprehensive study of its analytical properties]{dytso2018analytical}. 
The EP generalises the Gaussian ($q=2$) and Laplacian ($q=1$) densities \citep[e.g.][]{nadarajah2005generalized} and corresponds to the $q$-norm optimisation problem with bridge penalty $\text{pen} (\theta) =  - \left(c|\theta|\right)^q$ \citep{frank1993statistical}.
More preciesly, $\text{EP}(x \mid q=2,\kappa)$ is a Normal distribution with zero mean and variance $\frac{1}{\kappa} - 1$, while $\text{EP}(x \mid q=1,\kappa)$ is a Laplace distribution with zero mean and variance $2\left(\frac{1-\kappa}{\kappa}\right)^{2}$. 
While the EP is often referred to as a Generalised Normal distribution, we want to stress the potential of viewing its folded version as a Generalised Gamma distribution, as detailed below. 
Moreover, the smaller the shape parameter $q$, the flatter the tails of the EP distribution become.

The GB2 is also known as Inverted Beta distribution, Generalised Beta prime distribution, or scaled F-distribution.
It includes many common distributions as special cases, such as Gamma, Lognormal, Weibull, and Exponential~\citep{mcdonald1984some}.

The cumulative distribution function of a GG is given in terms of the lower incomplete gamma function $\gamma$ by
$
F_{\text{GG}}(x \mid a, d, q) = \frac{\gamma(d/q, (x/a)^q)}{\Gamma(d/q)}
$,
where $\gamma(d/q, (x/a)^q) = \int_{0}^{(x/a)^q}t^{d/q}e^{-t}\text{d} t$.

Note that the kernel of the GG is equivalent to the kernel of the EP when $d=1$, up to reparametrisation. This implies that the GG coincides with the folded version of an EP.

The additional scale parameter $\xi$ of the G3B distributions allows for a much wider variety of shapes than the standard Beta distribution (corresponding to $\xi = 1$). 
The G3B was introduced by \cite{libby1982multivariate}  \citep[see][Section 4.6, for applications, properties and generalisations, e.g., a five parameter version]{nadarajah2007multitude}.
It has been used in the Bayesian literature for estimating the ratio of two variances  \citep{gelfand1988estimation}. In the shrinkage literature, the G3B has been used as mixing density for SMNs by \cite{armagan2011generalized}.
The G3B is a special case of the Gauss Hypergeometric
and the Compound Confluent Hypergeometric distribution \citep[see, e.g.,][]{nadarajah2007multitude}.

\section{Gambel distributions}\label{sec:Gambel}

\subsection{Scale Beta mixtures of Exponential Power laws}
Scale mixtures of Exponential Power distributions generalise scale mixtures of Normals. They have been considered by~\cite{zhang2012ep} as shrinkage priors, with mixing density from the class of Generalised Inverse Gaussian densities, resulting in marginal priors that are nonsingular at zero. We consider such scale mixtures with three-parameter Beta densities \citep{libby1982multivariate}, resulting in marginals that can facilitate both singular and nonsingular behaviour at zero, and encompassing Horseshoe-like priors. 

In the following, we use the notation $ \Gamma \scriptsize{ \left[\begin{array}{c} a_1, a_2, \ldots, a_m \\ b_1, b_2, \ldots, b_n \end{array}
\right]} =
\frac{\Gamma(a_1)\Gamma(a_2) \cdots \Gamma(a_m)}{\Gamma(b_1)\Gamma(b_2) \cdots \Gamma(b_n)}
$.
We denote by $\Psi\left( x,y,z\right)$ the Tricomi confluent hypergeometric function, given for $y=1+\frac{1}{q}-b\notin \mathbb{Z}$ by
\begin{equation*}
	\Psi\left( x,y,z\right) = \frac{\pi}{\sin(\pi y)} \left\{\frac{M(x,y,z)}{\Gamma(1+x-y)\Gamma(y)} - 
	z^{1-y}\frac{M(1+x-y,2-y,z)}{\Gamma(x)\Gamma(2-y)}\right\}
\end{equation*}
with $M(x,y,z)$ denoting the Kummer confluent hypergeometric function \citep{silverman1972special}. Further,
we refer for all proofs to the Supplement. 
\begin{theorem}\label{th:density}
Given parameters $q,a,b,\xi>0$, consider a random variable $\theta$, taking values in $\mathbb{R}$, such that $\theta\mid \kappa \sim \text{EP}(\theta \mid q,\kappa)$ and $\kappa \sim \text{G3B}(\kappa \mid a,b,\xi)$. Then the marginal density of $\theta$ is
\begin{equation} \label{eq:Gambelden}
	f(\theta\mid q,a,b,\xi) = 
	\sftwo{\Gamma}{a+\tfrac{1}{q}}{\tfrac{1}{q}}{}
	\tfrac{q}{2\Beta(a,b)}\tfrac{\phi(q)}{\xi^{1/q}}
	\Psi\left( a+\tfrac{1}{q},1+\tfrac{1}{q}-b, \tfrac{(\phi(q)\,|\theta|)^q}{\xi}\right) \quad(\theta\in\mathbb{R})
\end{equation}

\end{theorem}
\begin{corollary}
    \label{th:centerandtailbehav}
The density $f(\theta\mid q,a,b,\xi)$ defined in Theorem~\ref{th:density} has a singularity at zero if and only if $b\leq 1/q$ and its tails behave like $|\theta|^{-aq-1}$.
\end{corollary}

By studying the Generalised Gamma, we will be able to relate the different hyperparameters to tail and origin behaviour. This has been pointed out by \cite{polson2012half,perez2017scaled,bai2018beta} for the Normal-mixture case, who remark in particular that singular and nonsingular behaviour can be facilitated by different choices of $a$. Additionally, $q$ can be used to steer the tail and origin behaviour.
Corollary \ref{th:centerandtailbehav} implies that the distribution tails follow a power law, whose behaviour is controlled by $a$ and $q$. Larger values of $a$ result in a faster decay to zero, making larger values of the regression coefficients less probable. 
Assigning a G3B for the shrinkage weight $\kappa$ is equivalent to specifying a GB2 on the shrinkage coefficient $\lambda^2$, as shown in the following proposition.

\begin{proposition}\label{prop:G3BtoGB2}
    Given parameters $q,a,b,\xi>0$, a random variable $\theta$ taking values in $\mathbb{R}$ is distributed according to
$\theta\mid \kappa \sim \text{EP}(\theta \mid q,\kappa)$, with $\kappa \sim \text{G3B}(\kappa \mid a,b,\xi)$ if and only if $\theta\mid \kappa \sim \text{EP}(\theta \mid q,\kappa)$, with 
$\lambda^2:= \frac{1-\kappa}{\kappa} \sim \text{GB2}(\lambda^2 \mid 1, \xi, b, a)$. 
\end{proposition}

A multitude of shrinkage priors can be obtained by varying the parameters $q,\xi,a$ and $b$. 
The original Horseshoe (HS) prior of \cite{carvalho2009handling,HorseshoeB}, defined by setting a half-Cauchy random variable as the standard deviation $\lambda_i$ of a Gaussian distribution (or equivalently, setting a Beta random variable as the local shrinkage weight), $\theta_i \mid \lambda_i, \tau \sim N(0, \lambda_i^2\tau^2)$ with 
$\lambda_i \sim C^{+}(0,1)$,
leads to a fat-tailed prior with unbounded density at the origin. Setting for simplicity $\tau = 1$, the HS prior can be reframed as a scale Beta mixture of exponential power laws as $\theta_i \mid \kappa_i\sim \text{EP}(\theta_i \mid 2,\kappa_i)$, with $\kappa_i \sim \text{G3B}(\kappa_i \mid 1/2, 1/2, 1)$.
Thus, our framework reduces to the original HS when $a = b = 1/q = 1/2$ and $\xi=1$.
\begin{figure}[t]
    \centering
\includegraphics[width=0.75\textwidth]{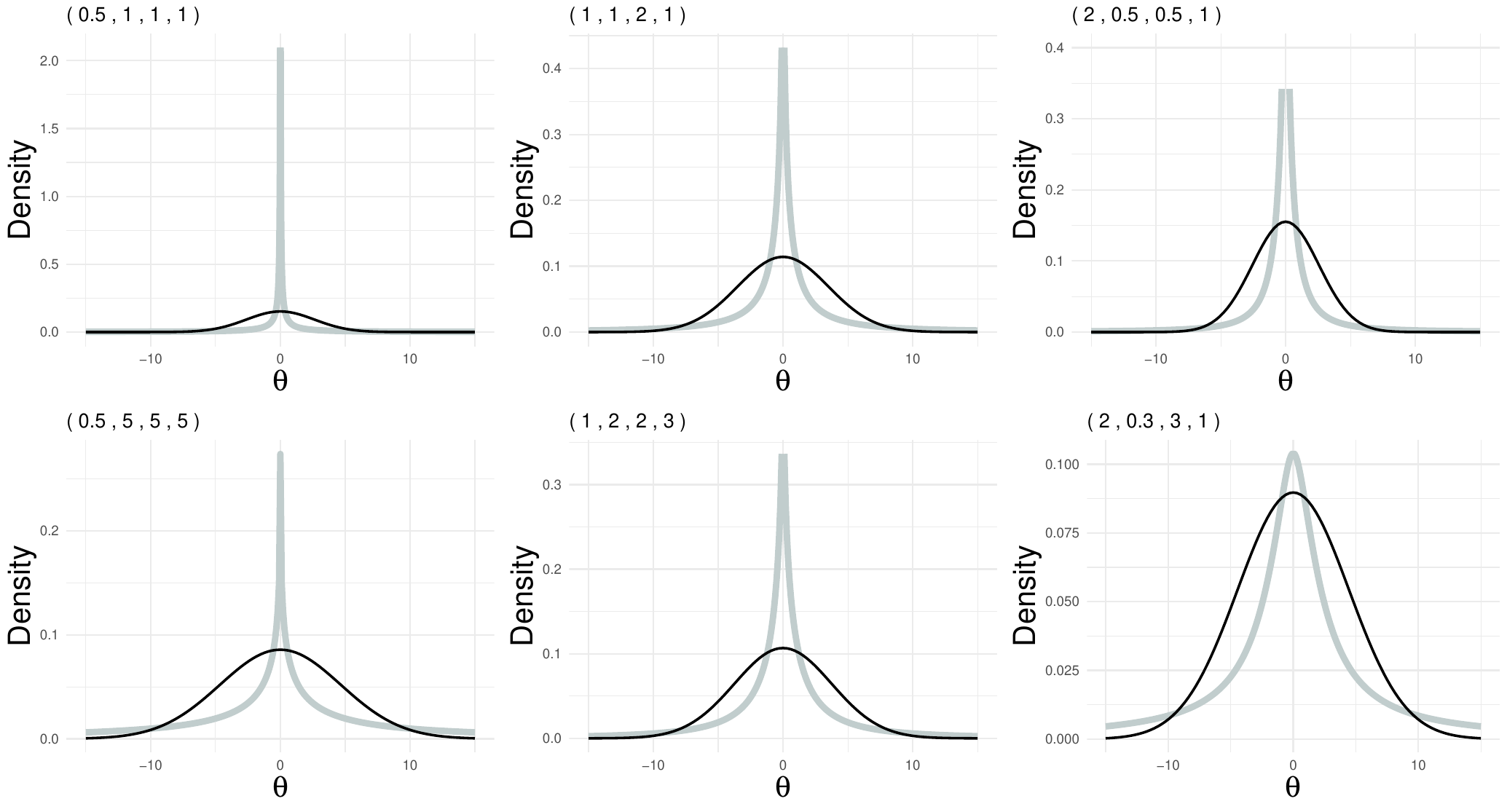}
    \caption{Density plots of the scale mixtures in Theorem \ref{th:density} for different hyperparameters (light solid line) versus the Normal density with same variance (dark thin line).}
    \label{fig:dens}
\end{figure}
Note that  \citet{zhang2012ep} propose the EP-GIG, which is also a sparsity inducing prior. This prior is a scale mixture distribution, where the parameter $\kappa/(1-\kappa)$ follows a Generalised Inverse Gaussian distribution, which includes the Gamma distribution as a special case. For this special case, the EP-GIG and the proposed scale mixture in \eqref{eq:Gambelden} coincide.

\subsection{Gamma-Beta quotient law}
In order to study the properties of the scale mixtures introduced in the previous section, it is convenient to focus on their folded versions on $\mathbb{R}^+$. All results derived for the folded versions can be easily translated in terms of the original shrinkage prior of interest. Nonetheless, in Theorem \ref{thm:ratio}, we show that the folded version on $\mathbb{R}^+$, to which we refer as Gambel distribution, admits a representation in terms of the ratio between two random variables, where the numerator and the denominator are distributed according to a generalised Gamma and a generalised Beta of the second kind, respectively. The cumulative distribution function and moments of the Gambel distribution are given in Propositions \ref{prop:cdf} and \ref{prop:moments} below. Formally, given $q,a,b,\xi>0$ and the density $f(\theta\mid q,a,b,\xi)$ defined in Theorem~\ref{th:density}, consider 
for $\theta\in\mathbb{R}^+$ the folded density 
(where, as above,  $\phi(q) =\sqrt{\Gamma(3/q) / \Gamma(1/q)}$)
\begin{align*}
\text{Gambel}(\theta\mid q,a,b,\xi) 
&:= 2 \times f(\theta\mid q,a,b,\xi) \\
&=\sftwo{\Gamma}{a+\tfrac{1}{q}}{\tfrac{1}{q}}{}
	\tfrac{q}{\Beta(a,b)}\tfrac{\phi(q)}{\xi^{1/q}}
	\Psi\left( a+\tfrac{1}{q},1+\tfrac{1}{q}-b, \tfrac{(\phi(q)\,\theta)^q}{\xi}\right)
\end{align*}
\begin{theorem}\label{thm:ratio}
If $X\sim \text{GG}(\xi^{1/q},d=1,q)$ and $Y\sim GB2(q,\phi(q),a,b)$, with $X \perp Y$, then 
\[
\theta = X/Y\sim \text{Gambel}(\theta \mid q,a,b,\xi)
\]
\end{theorem}

\begin{proposition}\label{prop:cdf}
The cumulative distribution function of the $\text{Gambel}$ distribution is 
\begin{multline}\label{cdf}
F(\theta; q,a,b,\xi) =\left(\tfrac{\phi(q)\,\theta}{\xi^{1/q}}\right)
\tfrac{  q\,\Beta\left(a+\tfrac{1}{q}, b-\tfrac{1}{q}\right) }{\Gamma\left(\tfrac{1}{q}\right)\Beta(a,b)}
\sffour{\, _2F_2}{\tfrac{1}{q},}{ a +\tfrac{1}{q} ;}{ 1+\tfrac{1}{q},}{ 1+\tfrac{1}{q} -a;}{ \left(\tfrac{\phi(q)\,\theta}{\xi^{1/q}}\right)^q  }  \\
+\left(\tfrac{\phi(q)\,\theta}{\xi^{1/q}}\right)^{qb}\tfrac{\Gamma\left(\tfrac{1}{q} - b\right) }{b \Gamma\left(\tfrac{1}{q}\right)\Beta({b},{a})} 
\sffour{\,  _2F_2}{b,}{ b+a;}{ 1+b,}{ 1+b-\tfrac{1}{q};}{\left(\tfrac{\phi(q)\,\theta}{\xi^{1/q}}\right)^q }
\quad(\theta>0)
\end{multline}
where $_2F_2$ is the generalised hypergeometric function.
\end{proposition}

\begin{proposition}\label{prop:moments} 
If $\theta\sim \text{Gambel}(\theta\mid q,a,b,\xi)$, the $k$-th moment  of $\theta$ exists whenever $k < aq$ and, in this case, is 
\begin{align}\label{eq:moments}
\mathbb{E}(\theta^k)& =   \left(\frac{\xi^{1/q}}{\phi(q)}\right)^k
\sfsix{\Gamma}{\tfrac{1+k}{q},}{a-\tfrac{k}{q},}{b+\tfrac{k}{q}}  
{\tfrac{1}{q},}{a,}{b}{}
\end{align}
In particular, if $\theta$ is distributed according to the scale mixture in Theorem~\ref{th:density} (i.e., the unfolded Gambel, $\theta\mid \kappa \sim \text{EP}(\theta \mid q,\kappa)$ with $\kappa \sim \text{G3B}(\kappa \mid a,b,\xi)$), the odd moments are zero whenever they exist. E.g.\ the mean is zero whenever $a>1/q$. Moreover, the variance is 
$\xi^{2/q}
\sffour{\Gamma}{a-\tfrac{2}{q},}{b+\tfrac{2}{q}}{a,}{b}{}$
if $a>2/q$; the fourth moment is
$\frac{\xi^{4/q}}{\phi(q)^4}
\sfsix{\Gamma}{\tfrac{5}{q},}{a-\tfrac{4}{q},}{b+\tfrac{4}{q},}{\tfrac{1}{q},}{a,}{b}{}$ if $a>4/q$;
and the kurtosis, in this case, is
$\frac{\Gamma(1/q)\Gamma(5/q)}{\Gamma^2(3/q)}\frac{\Gamma(a-4/q)\Gamma(b+4/q)}{\Gamma^2(a-2/q)\Gamma^2(b+2/q)}\Gamma(a)\Gamma(b)$,
which can be approximated as
$\left(\frac{3a(qb+2)}{b(qa-4)}\right)^{2/q}$.
%\MDI{Chosen $a,b$ by checking the kurtosis.}
\end{proposition}

\begin{remark}
\cite{nadarajah2012gamma} study the properties of the Gamma Beta Ratio distribution, a special case of the Gambel distribution in Theorem~\ref{thm:ratio}.
If $X\sim\textrm{Gamma}(\beta,\lambda)$ and $Y\sim\textrm{Beta}(a,b)$, then $Z=X/Y$
has density and cdf
\begin{align}\label{KotzGBratio}
f_Z(z)=\frac{\lambda^{\beta}\Beta(\beta+a,b)}{\Gamma(\beta)\Beta(a,b)}
z^{\beta-1}\,_1F_1(\beta+a;\beta+a+b;-\lambda z)
\quad(z>0)
\end{align}
\[
F_Z(z) =\frac{\Beta(b,a+\beta)(\lambda z)^{\beta}}{\Gamma(\beta+1)\Beta(a,b)}\,_2F_2(\beta,a+\beta;\beta+1,a+b+\beta;-\lambda z)
\]
 In \eqref{KotzGBratio}, $\beta=1$ reduces the Gamma to the Laplace distribution and we obtain a scale-mixture of a Laplace with a beta-prior on the reciprocal of the standard deviation. This leads to a distribution with long fat tails.   
\end{remark}

\section{Characterization of scale mixture priors}\label{sec:geometry}
\subsection{Decreasing hazard and concave survival function}
Given the survival function 
$
S(x):=
\mathbb{P}(X> x) = 1 -F(x)
$, the \emph{hazard function} is the non-decreasing function $R(x):=-\ln S(x)$ and the corresponding \emph{hazard rate} is
\[
r(x) := R'(x) = -\frac{S'(x)}{S(x)}= \frac{f(x)}{S(x)} = \frac{\mathbb{P}(X\in \text{d} x)}{\mathbb{P}(X\geq x)}
\] 
The hazard rate, also called \emph{failure rate} in reliability, \emph{mortality of claims} in actuarial applications and \emph{intensity function} in extreme value theory, can be used to distinguish among the various skewed fat right-tail distributions. The lower the hazard rate, the more skewed and ``dangerous'' the distribution is \citep{kleiber2003statistical}. 
Many distributions in reliability applications have increasing hazard rates, at least from some argument onwards. 
If $X$ is the time to a certain event and the hazard rate is increasing then the conditional probability of experiencing the event given survival up to time $x$  increases with time (e.g. elderly mortality). 
A constant hazard function means the chances of experiencing the event in case of survival up to $x$ do not change with $x$, i.e. the failure risk stays the same (e.g. adult mortality). Contrary, when the hazard is decreasing, the probability of experiencing the event decreases with time (e.g. child mortality).

Since
$r'(x)=%\tfrac{S(x) f'(x) + f^2(x)}{S^2(x)} = \frac{f'(x)}{S(x)} + \frac{f^2(x)}{S^2(x)} =
 \frac{f'(x)}{S(x)} + r^2(x)
$,
the monotonicity of the hazard rate is characterized by $\pen'(x) = (-\ln f(x) )'$. We have that $r(x)$ is decreasing (resp.\ constant, increasing) if and only if $\pen'(x)$ is larger than (resp.\ equal to, smaller than) $r(x)$.

In the reliability literature \citep[e.g.][]{mcdonald1988hazard} the notation $\eta(x):=(-\ln f(x) )' $ is used, i.e.\ $\eta=\pen '$. The hazard rate is increasing (decreasing) whenever $\eta'(x)>0$ ($<0$) for all $x>0$. (\cite{glaser1980bathtub} derives this from noticing $(1/r)'=(1/r)\eta - 1$. In case $\eta$ has one simple zero, the situation splits up in further (non-monotonic) cases.
\begin{figure}
    \centering
    \includegraphics[width=0.75\linewidth]{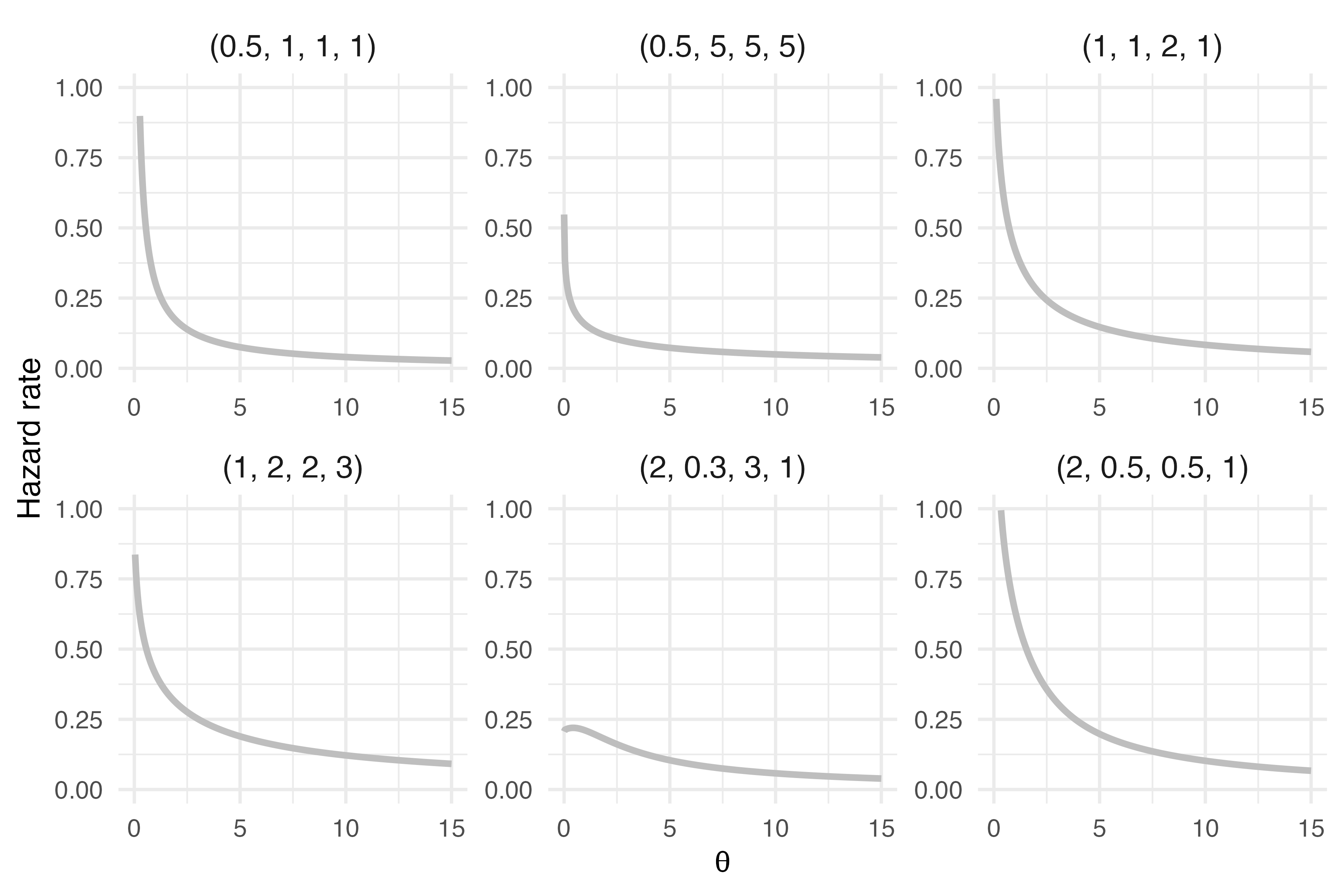}
    \caption{Hazard rates of the folded Gambel distribution for various hyperparameters.}
    \label{fig:hGambel}
\end{figure}
The monotonicity of the hazard rate is also determined by the second derivative of the hazard function $R$. In particular, $r$ is decreasing (resp.\ increasing) in $x$ if and only if $R$ is concave (resp.\ convex) in $x$. By definition of $R$, this means in terms of the tail function that
$r(x)$ is decreasing (resp.\ increasing) in $x$ if and only if
$\ln (1-F(x))$ is convex (resp.\ concave) in  $x$.

\begin{theorem}\label{thm:hazard}
Consider the hazard rate $r(x)$ of a Gambel distribution with parameters $q,a,b,\xi >0$. 
\begin{enumerate}
\item If $q<1$, the hazard rate $r(x)$ is decreasing on $(0,\infty)$.
\item If $q\geq 1$, then there exists some $x_0$ such that $r(x)$ is decreasing on $(x_0,\infty)$.
\end{enumerate}
\end{theorem}

Figure \ref{fig:hGambel} shows hazard rate curves for the folded Gambel distribution for several choices of hyperparameters. In our context, the interpretation of a decreasing hazard rate is that a regression coefficient that is not shrunk to zero can assume large values. 

\subsection{Lorenz ordering}
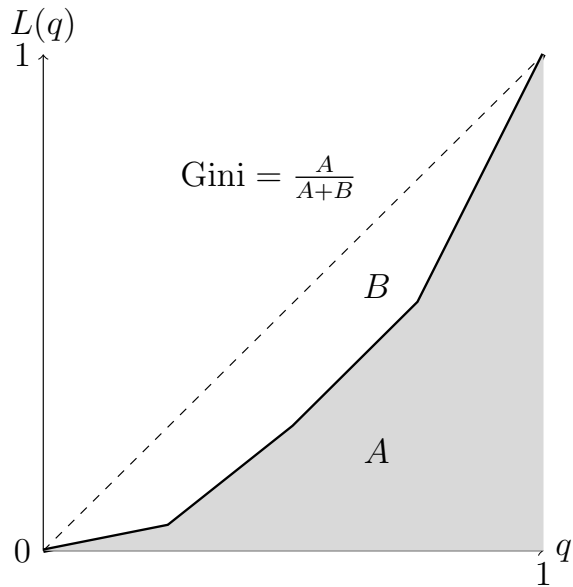
\begin{figure}[tb!]
  \centering
  \resizebox{0.45\textwidth}{!}{
  \begin{tikzpicture}
    % Axes
    \draw[->] (0,0) -- (6,0) node[right] {$q$};
    \draw[->] (0,0) -- (0,6) node[above] {$L(q)$};
    
    % Perfect equality line
    \draw[dashed] (0,0) -- (6,6);
    
    % Lorenz curve
    \draw[ultra thick] (0,0) -- (1.5,0.3) -- (3,1.5) -- (4.5,3) -- (6,6);
    
    % Shading under the Lorenz curve
    \fill[gray!30] (0,0) -- (1.5,0.3) -- (3,1.5) -- (4.5,3) -- (6,6) -- (6,0) -- cycle;
    
    % Gini coefficient
    \draw (1.5,4.5) node[right] {$\text{Gini} = \frac{A}{A+B}$};
    
    % Labels
    \draw (6,0) node[below] {1};
    \draw (0,0) node[left] {0};
    \draw (0,6) node[left] {1};

    % Area labels
    \draw (4,1.2) node {$A$};
    \draw (4,3.2) node {$B$};
  \end{tikzpicture}
  }
  \caption{\label{fig:lorenz-example}Illustration of Lorenz Curve and Gini Index:
The thick line depicts a Lorenz curve.
The Gini index is calculated as the ratio of the area of the shaded gray region ($A$) to the total area under the perfect equality line ($A+B$).}
\end{figure}

Lorenz curves \citep{lorenz1905methods} are functional measures and graphical representations widely employed in economics for the study of wealth distributions, where they  provide valuable insight into wealth inequality. They measure the concentration of income or wealth across a population and allow to visually compare it to the situation of perfect equality, in which each segment of the population holds an equal share of the total income or wealth.  Given an observed variable $Y$ corresponding, for instance, to a measure of wealth, the sample Lorenz curve of a sample $(Y_1,\ldots,Y_n)$ is computed by calculating the relative cumulative wealth $C(x) = (\sum_{i=1}^n Y_{i}\mathbbm{1}_{(Y_{i}\leq x)})/(\sum_{i=1}^n Y_{i})$, i.e., the percentage of wealth detained by all subjects that have at most $x$, and plotting it against the empirical cumulative distribution function $F(x) = \tfrac{1}{n}\sum_{i=1}^n\mathbbm{1}_{(Y_{i}\leq x)}$. Hence, the Lorenz curve is $L(p) = C(F^{-1}(p))$, where $F^{-1}(p)=\inf\{x:F(x)\geq p\}$ is the generalised inverse of the cumulative distribution. See~Figure~\ref{fig:lorenz-example} for a graphical example. 
Given an absolutely continuous distribution with density $f(x)$, the Lorenz curve of $f(x)$ can be defined as 
$L(p) = {\int\limits_{0}^{F^{-1}(p)}x f(x) \text{d} x}
\,/{\int\limits_{0}^{+\infty}x f(x) \text{d} x}$,
while the Gini coefficient is given by
\begin{equation}\label{eq:gini_coeff}
G = \frac{\int_{0}^{+\infty}\int_{0}^{+\infty}f(x)f(y)|x-y|\,\text{d}x \,\text{d}y}{2\int_{0}^{+\infty}x f(x) \,\text{d} x}
\end{equation}

When comparing Lorenz curves, a partial ordering called \emph{Lorenz ordering} or \emph{Lorenz dominance} is established. Lorenz ordering provides a formal way to compare concentration levels across different distributions.
Two distributions functions $F_1$ and $F_2$ with corresponding Lorenz curves $L_1$ and $L_2$, respectively, are said to exhibit Lorenz ordering if either $L_1(p) \leq L_2(p)$ or $L_1(p) \geq L_2(p)$, for all $p \in [0,1]$.
In the latter case $L_1$  is said to Lorenz dominate $L_2$, implying that its distribution exhibits less inequality than the distribution represented by $L_2$.
Since Lorenz ordering establishes only a partial order, certain distributions can be incomparable in terms of the inequality levels described by Lorenz curves.

In the context of reliability theory, numerous similar concepts have been developed to compare the time-to-failure patterns of different components or systems. These concepts establish a partial ordering among the probability distributions associated with them. Interestingly, a notable resemblance exists between these partial orders in reliability theory and the ones of Lorenz ordering \citep{chandra1981relationships, kochar2006lorenz}.

In general, Lorenz curves, Lorenz orderings, and Gini indexes measure concentrations of distributions, which in the context of shrinkage priors can be used as an indicator of the level of shrinkage. Taking as reference a uniform improper distribution over the real line for each of the coefficients in the regression, $p(\theta)\propto 1$ results in no penalization of the estimates and thus in the MLE coefficients. The larger the ``distance'' from such prior, the higher the shrinkage. Table \ref{table:Gini} shows the Gini coefficient of the folded Gambel distribution for several hyperparameter values. As a closed-form expression for the Gini coefficient is not attainable, we employed numerical integration of \eqref{eq:gini_coeff}.
\begin{table}[h]
\centering
\caption{Gini indexes for the Gambel distribution with the various  hyperparameter choices depicted in Figures~\ref{fig:dens} and \ref{fig:lorenz}.}
%\resizebox{0.4\textwidth}{!}{
\begin{tabular}{ccccc}
\toprule
$q$ & $a$ & $b$ & $\xi$ & Gini index  \\ \midrule
0.5 & 1   & 1   & 1          & 0.002 \\
1   & 1   & 2   & 1          & 0.077 \\
2   & 0.5 & 0.5 & 1          & 0.120 \\
0.5 & 5   & 5   & 5          & 0.235 \\
1   & 2   & 2   & 3          & 0.364 \\
2   & 0.3 & 3   & 1          & 0.002\\
\bottomrule
\end{tabular}%}
\label{table:Gini}
\end{table}

\begin{figure}[htb]
    \centering    \includegraphics[width=0.75\textwidth]{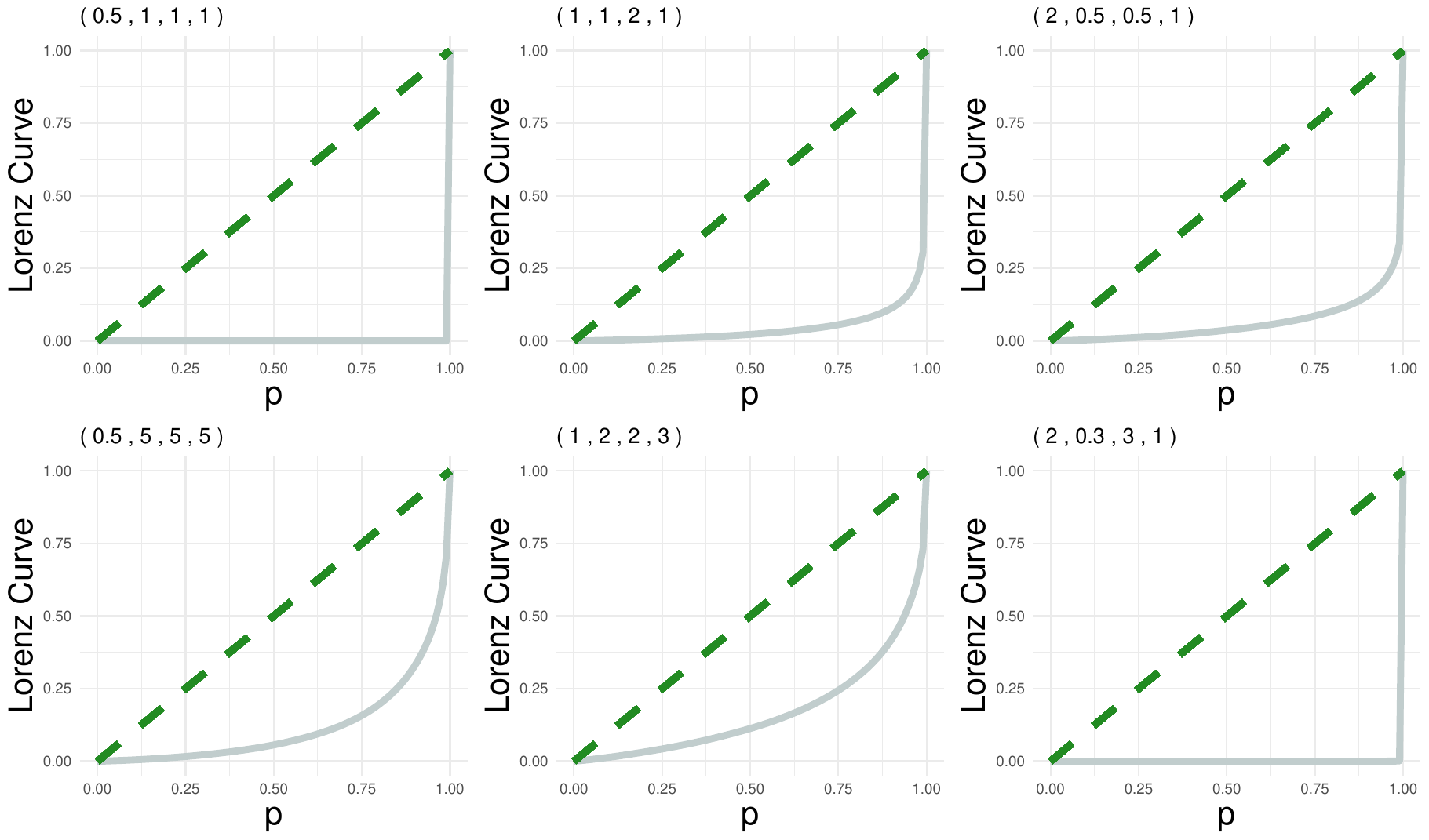}
    \caption{Lorenz curves of scale mixtures for different hyperparameters. Dashed lines denote the improper prior distribution of perfect equality.}
    \label{fig:lorenz}
\end{figure}

\begin{theorem}\label{th:Lorenz}
Defining $\theta^{\star}:= F^{-1}(p)$, the Lorenz curve of a Gambel distribution with parameters $q,a,b,\xi>0$ is 
\begin{scriptsize}
\begin{align*}
   L(p) &= C\Bigg(\theta^* \xi^b (bq+1)\phi(q)\Gamma\!\left(a+\frac{1}{q}\right) 
   \Gamma\!\left(b-\frac{1}{q}+1\right) {_2F_2}\!\left(a+\frac{1}{q},\frac{2}{q};-b+\frac{1}{q}+1,1+\frac{2}{q};
   \frac{\phi(q)^q\!(\theta^*)^q}{\xi}\right) \\
   &\quad -2 \xi^{\frac{1}{q}} \Gamma(a+b)(\theta^*)^{bq} 
   \Gamma\!\left(-b+\frac{1}{q}+1\right)\phi(q)^{bq} {_2F_2}\!\left(a+b,b+\frac{1}{q};b-\frac{1}{q}+1,b+\frac{1}{q}+1;
   \frac{\phi(q)^q\!(\theta^*)^q}{\xi}\right)\!\Bigg) \\
   &\text{where }C = \frac{\pi\theta^* q \sqrt{\Gamma\!\left(\frac{3}{q}\right)} 
   \xi^{-\frac{bq+2}{q}}\csc\!\left(\pi\!\left(b-\frac{1}{q}\right)\right)}
   {2(bq+1)\Gamma\!\left(\frac{1}{q}\right)^{\!3/2} \Gamma(a+b)B(a,b)}\times \frac{1}{\Gamma\!\left(b-\frac{1}{q}+1\right) 
   \Gamma\!\left(-b+\frac{1}{q}+1\right)} \\
   &\quad\times \frac{1}{_3F_3\!\left(a-\frac{1}{q},b+\frac{1}{q},\frac{2}{q};
   a,b,\frac{1}{q};1\right)}
\end{align*}
\end{scriptsize}
\end{theorem}

\section{Algorithms}\label{sec:algorithms}

In this section we describe two Markov chain Monte Carlo (MCMC) algorithms for posterior inference. We consider separately two cases: the low-dimensional case $p<n$, and the high-dimensional case $p\geq n$. For the Horseshoe prior, more efficient algorithms can be found, for instance, in \citet[][]{bhattacharya2021mcmc}.

\subsection{MCMC for the case $p<n$} \label{sec:MCMC1}

In case $p<n$, our proposed algorithm is similar to the Bayesian Bridge algorithm \citep{polson2012local}, which it can be easily extended to include as a particular case. Before presenting the hierarchical formulation that serves as the basis for our MCMC algorithm, we discuss approaches taken by previous authors. \cite{mallick2014new} use a mixture of uniforms \citep[an idea from][]{boris2008scale} to represent the EP. The EP can be expressed as SMN
$e^{-|z|^q}\propto\int_0^{\infty}\tfrac{1}{\sqrt{2\pi}s}e^{-z^2/(2s)}g_{q/2}\left(\tfrac{1}{2s}\right)\,\text{d}s$
\cite[see the seminal work on the Bayesian Lasso by][]{park2008bayesian}, where $g_{q/2}$ is the density of a positive stable random variable with index $q/2$  which in general has no closed form expression \cite[see][]{west1987scale}. The resulting conditional distributions are closely related to certain exponential dispersion models, but it is unclear whether an efficient Gibbs sampler can be developed. Moreover, the resulting conditional posterior and prior are not known in closed form \citep{polson2014bayesian}. \cite{polson2019bayesian} therefore show how to represent the GG distribution as a scale mixture of Bartlett-Fejer kernels, and this allows building a Gibbs sampler (though only for the case $p\leq n$ and full rank $\bX$).
Additionally, \cite{alhamzawi2018new} follow \cite{mallick2014new} almost verbatim, but write the Laplace density as a scale mixture of truncated Normal distributions with exponential mixing densities. They claim that this leads to better, more stable and faster convergence of the Gibbs sampler in case $p\gg n$ or in the presence of multicollinearity, i.e.\ when $\bX'\bX$ is singular or nearly so. Further literature includes \cite{NardonPianca} and  \cite{kalke2013simulation}, who also built an R package \emph{pgnorm} for simulation. We next provide a result that is useful for devising an MCMC scheme. It is an extension of Proposition 1 in \citet{armagan2011generalized}. 
\begin{proposition}\label{Thm:fullconditionals}
     Let $\theta\mid\kappa \sim \EP(q,\kappa)$, $\kappa \sim \TPB(a,b,\xi)$, and $\displaystyle \gamma = \left(\tfrac{1-\kappa}{\kappa}\right)^{1/q}$.

\begin{itemize}
\item[(i)]  Then $	\gamma\mid\lambda \sim \GG(\lambda^{-1/q},qa,q)$ and $\lambda  \sim\mbox{Gamma} (b,\xi)$. 

\item[(ii)]   The marginal distribution of $\gamma$ is the generalized Beta-prime distribution 
\[ 
\pi(\gamma) = q\xi^b\frac{\Gamma(a+b)}{\Gamma(a)\Gamma(b)}\gamma^{qa-1}[\xi+\gamma^q]^{-(a+b)}
\]

\item[(iii)] The conditional posterior distribution of $\gamma\mid \theta$ is 
\[
\pi (\gamma \mid \theta) \propto  \exp\left\{
 - \frac{\phi(q)^q|\theta|^q}{\gamma^q} \right\} \gamma^{qa-2}  \frac{1}{[\xi+\gamma^q]^{a+b}}%,
\] 
\end{itemize}    
\end{proposition}

Our algorithm is an extension of the work by 
\cite{mallick2014new}, who expressed the Laplace density as a scale mixture of uniform distributions (see also \cite{boris2008scale}). This representation of the Laplace density is given in the following proposition (proof given in the Supplement).
\begin{proposition}\label{Thm:GGDasSMU}
The EP distribution with parameters $q,\kappa$ can be obtained as a scale mixture of a uniform distribution, with a Gamma distribution mixing density, i.e.\
$$
\pi(\theta)=
\int\limits_{-\gamma u^{1/q}<\theta<\gamma u^{1/q}} \frac{1}{2\gamma u^{1/q}}f(u\mid 1+1/q,\phi(q)^q)\,\textup{d}u.
$$
\end{proposition}

The EP thus has the  hierarchical representation $\theta \mid u \sim \mbox{Uniform}(-\gamma u^{1/q},\gamma u^{1/q})$, with 
$u \sim \mbox{Gamma}(1+1/q,\phi(q)^q)$, and a such our hierarchical model is:
\begin{align*}
\by  \mid \bX,\bbeta,\sigma^2 & \sim \mbox{N}(\bX\bbeta,\sigma^2\bI_n)\\
\sigma^2 &\sim \pi(\sigma^2) \quad\text{ (usually $\mbox{IG}(\alpha,\theta)$)}\\
\theta_j  \mid u_j,\lambda^2_j &\sim \mbox{Uniform}(-\gamma_j u_j^{1/q},\gamma_j u_j^{1/q})\\ 
u_j &\sim \mbox{Gamma}(1+1/q,\phi(q)^q)\\
%\bbeta  \mid \bu,\sigma^2 &\sim \prod_{j=1}^p  \textrm{Uniform}(-\sqrt{\sigma^2}u_j,\sqrt{\sigma^2}u_j)\\
	\gamma_j\mid\lambda_j \sim \mbox{GG}(\lambda_j^{-1/q},&qa,q), \quad
\lambda_j \sim\mbox{Gamma} (b,\xi)   
\end{align*}

The Gibbs sampler is built by sampling from the full conditionals (shown in Supplement) of each random variable in the previous hierarchical model. 

\subsection{MCMC for the case $p \geq n$} \label{sec:MCMC2}
When $p \geq n$, the strategy is to write the EP  as a scale mixture of normals so that, in case of singular $\bX'\bX$, we can employ an MCMC algorithm in the spirit of \cite{alhamzawi2018new}. \cite{polson2012local} show that %[Theorem 1] 
\begin{align}\label{PolsonScott}
k\exp\{-\nu|\theta|^q\}
&\propto
\int_0^{\infty} 
\textrm{N}(\theta\mid 0,T^{-1}) T^{-1/2} \tfrac{1}{\nu^{2/q}}
p_{q/2}\left(\tfrac{T}{\nu^{2/q}}\right) 
\textrm{d}T\\
&=
\tfrac{1}{\sqrt{2\pi}}\int_0^{\infty}
\exp\left\{-\tfrac{T\theta^2}{2}\right\}  \tfrac{1}{\nu^{2/q}}
p_{q/2}\left(\tfrac{T}{\nu^{2/q}}\right)
\textrm{d}T   \label{PolsonScott2}
\end{align}
where $k$ is the normalising constant of the EP density, and $p_{q/2}$ is the standardized and positive $\tfrac{q}{2}$-stable density (here \emph{positive} means skewness parameter one, so that the distribution is supported on the positive real line). 
In our context, we get 
\begin{align*}
\pi_{q}(\theta | \gamma) &= 
\frac{q\phi(q)}{2\gamma\Gamma(1/q)} \exp\left\{ - \left( \frac{\phi(q)|\theta|}{\gamma}\right)^q\right\}\\
&=k
\int_0^{\infty}
\textrm{N}(\theta\mid 0,T^{-1}) T^{-1/2} \frac{\gamma^2}{\phi(q)^2}
p_{q/2}\left(\frac{T\gamma^2}{\phi(q)^2}\right)
\textrm{d}T   \\
&= k
\int_0^{\infty}
\textrm{N}\left(\theta\mid 0,\frac{\gamma^2}{\phi(q)^2\tau}\right) \frac{\gamma}{\phi(q)\tau^{1/2}}
p_{q/2}\left(\tau\right)
\textrm{d}\tau  \quad[\text{substitute } \tau:=T\gamma^2/\phi(q)^2] 
\end{align*}
which we use to obtain the following hierarchical model:
\begin{align*}
\by  \mid \bX,\bbeta,\sigma^2 & \sim \mbox{N}(\bX\bbeta,\sigma^2\bI_n)\\
\sigma^2 &\sim \pi(\sigma^2) \quad\text{ (usually $\mbox{IG}(\alpha,\theta)$)}\\
\theta_j  \mid \gamma_j, \tau_j &\sim \mbox{N}\left(0,\gamma_j^2/(\phi(q)^2\tau_j)\right)\\ 
p(\tau_j) & \sim \tau_j^{-1/2} p_{q/2}(\tau_j) \quad(q/2<1) \\
	\gamma_j\mid\lambda_j \sim \mbox{GG}(\lambda_j^{-1/q},qa,&q) \quad \lambda_j \sim\mbox{Gamma} (b,\xi)   
\end{align*}

The Gibbs sampler is built by sampling from the full conditionals (shown in Supplement) of each random variable in the previous hierarchical model. 

\begin{example}
For the lasso ($q=1$), the mixture \eqref{PolsonScott2} is with the $\tfrac{1}{2}$-stable distribution, which has density  
$
p_{1/2}(t) = \frac{1}{\sqrt{2\pi t^3}} \exp\left\{-\frac{1}{2t}\right\}$ for $t>0$
\citep[see][VI.2]{feller2008introduction}.
Therefore, by \eqref{PolsonScott2}, \cite[same as][p. 293]{polson2012local} we get
\begin{align*}
k\exp\{-\nu|\theta|\}
&\propto
\int_0^{\infty}
\exp\left\{-\tfrac{T\theta^2}{2}\right\} 
\tfrac{\nu}{\sqrt{2\pi T^3}} \exp\left\{-\tfrac{\nu^2}{2T}\right\}  
\,\textup{d}T 
\end{align*}    
\end{example}

\section{Consistency}\label{sec:asymptotic}
Consistency results are concerned with the behaviour of regression methods for growing number of predictors $p$ and sample size $n$. Denoting the dependence of the relevant notions by a subscript $n$, 
the posterior distribution of $\bbeta_n$ under a given prior with density $\pi(\bbeta_n)$ is called \emph{strongly consistent} if 
$\lim_{n\to\infty}\Prob(\|\bbeta_n - \bbeta_n^*\|>\epsilon \mid \by) = 0$ pr$\prescript{}{\bbeta^*_n}{}$-almost surely for any  $\epsilon>0$.

For scale mixtures of Normals, sufficient conditions for strong consistency have been shown by
\cite{armagan2013posterior} in the low-dimensional setting ($p\leq  n$), and have been extended to the high-dimensional setting ($p > n$) by \cite{bai2018high}. 
In this section, we extend strong consistency results to the entire class of scale Beta mixtures of Exponential Power laws introduced in Section~\ref{sec:Gambel}.

Following \cite{armagan2013posterior} and \cite{bai2018high}, we rely on the following mild and common assumptions.
We denote by $S\subseteq \{1,\ldots,p\}$ the support of the unknown
true vector $\bbeta^*$, i.e., $\beta_j^*>0$ if and only if $j \in S$, and by $s$ the cardinality of $S$.
Note that, besides $\bbeta^*$ and the parameter vector $\bbeta$, also the number of regressors $p$, and the sparsity level $s$ are understood as sequences depending on $n$, and denoted with an additional subscript.
The assumptions for the low-dimensional setting are
\begin{itemize}
\item[(L1)] $p_n\leq n$ $\forall n$  and $p_n=o(n)$ %\todoB{\cite{armagan2011generalized} does not impose  $p_n\leq n$,\cite{bai2018high} does}
\item[(L2)] $0<\Lambda_{\min}^2:=\liminf_{n\to\infty}\lambda_{\min}(\tfrac{1}{n}\bX_n'\bX_n)\leq
 \limsup_{n\to\infty}\lambda_{\max}(\tfrac{1}{n}\bX_n'\bX_n)=:\Lambda_{\max}^2<\infty$, where 
 $\lambda_{\min}$ and $\lambda_{\max}$ denote the minimum and maximum eigenvalues (in particular $\bX_n$ has full rank, i.e., $\bX_n'\bX_n$ is positive definite) 
\item[(L3)] $s_n = o\left(\tfrac{n}{\log n}\right)$ (i.e., $\lim_{n\to\infty} s_n(\log n)/n=0$) 
\end{itemize}

The assumptions for the high-dimensional setting are
\begin{itemize}
\item[(H1)]  $p_n>n$ $\forall n$ and $\ln(p_n)=O(n^d)$ for some $d\in(0,1)$
\item[(H2)] $\forall J \subset \{1,\ldots, p_n\} $ with $|J|\leq n: 0< \Lambda_{J,\min}^2:=\liminf_{n\to\infty}\lambda_{\min}(\tfrac{1}{n}\bX_J'\bX_J)$ and $
 \Lambda_{\max}^2<\infty$, where $\bX_J$ is the submatrix of the $J$-indexed columns of $\bX$.
\item[(H3)] $s_n=  o\left(\tfrac{n}{\log p_n}\right)$
\end{itemize}

As final condition (L4) and (H4), for both settings, the sizes of the coefficients are assumed to be uniformly bounded above in absolute value, i.e., $\sup_{n\in\N}\sup_{j\in S}|\theta^*_{n,j}| <\infty$.

\cite{armagan2013posterior} and \cite{bai2018high} show that sufficient for strong consistency is that the prior captures $\bbeta^*_n$ inside a ball of size depending on the eigenstructure of $\bX_n$ with sufficiently high probability. Specifically, given (L1)-(L2) \citep[see][Thm 1]{armagan2013posterior} or  (H1)-(H3)  \citep[see][Thm 2]{bai2018high} for the respective settings, sufficient for strong consistency is that 
\begin{align}\label{StronConsistencySuffCond}
\Prob\left(\|\bbeta_n-\bbeta_n^*\|<\tfrac{\Delta}{n^{\rho/2}}\right)
>\exp(-dn)
\end{align}
for all strictly positive $\Delta$ and $d$ up to sizes depending on $\epsilon, \Lambda_{\max}^2$ and $\Lambda_{\min}^2$ (resp. $\Lambda_{J,\min}^2$),  
and some $\rho>0$. 
To verify \eqref{StronConsistencySuffCond} for a sparsity promoting prior at hand, \cite{armagan2013posterior} start with the estimate
\begin{align}
\Prob\left(\|\bbeta_n - \bbeta_n^*\|<c_n\right) 
&=  \Prob\left(\sum_{j\notin S}\theta_{nj}^2  +  \sum_{j\in S}(\theta_{nj}-\theta^*_{nj})^2 <c_n^2\right)  \nonumber \\
&\geq   \Prob\left(\sum_{j\notin S}\theta_{nj}^2   <  \tfrac{(p_n-s_n)c_n^2}{p_n}\right)
\prod_{j\in S}  \Prob\left( |\theta_{nj}-\theta^*_{nj}| < \tfrac{c_n}{\sqrt{p_n}} \right)   \nonumber \\
&\geq   \left(1-\tfrac{p_n}{(p_n-s_n)c_n^2}\mathbb{E}\left(\sum_{j\notin S} \theta_{nj}^2\right)\right)
\prod_{j\in S}  \Prob\left( |\theta_{nj}-\theta^*_{nj}| < \tfrac{c_n}{\sqrt{p_n}} \right)   \nonumber
%\quad[\text{Markov}] 
\\
&\geq   \left(1-\tfrac{p_n}{c_n^2}\mathbb{E}(\theta_{nj}^2)\right)
\left(\tfrac{c_n}{\sqrt{p_n}} \pi\left(\sup_{j\in S}|\theta_{nj}^*| + \tfrac{c_n}{\sqrt{p_n}} \right)\right)^{s_n}  \label{eq:consistency_eq}
\end{align}
where the last inequality uses that $\pi$ is assumed symmetric, and decreasing on $\R^+$. Indeed, in this case, if $X\sim\pi$ and $r,m>0$, then 
$\Prob (|X-m|<r)\geq r \cdot \pi(m+r)$; the area of the rectangle $[m,m+r]\times[0,\pi(m+r)]$. In particular,
the range of the Gambel parameter $a$ in the following consistency result is to guarantee a finite second moment in \eqref{eq:consistency_eq}. 
\begin{theorem}\label{th:consistency}
Assuming (L1)-(L4) in the low-dimensional setting, and (H1)-(H4) in the high-dimensional setting, the Beta mixture of Exponential Power laws prior in Theorem \ref{th:density} with $a\in(2/q,\infty)$, $b\in(0,\infty)$ yields a strongly consistent posterior if $\xi= (C/(pn^\rho \log n))^{q/2}$ for some $C>0$.
\end{theorem}

\section{Simulation study}\label{sec:sim_study}

In this section, we conduct two simulation studies where we perform Bayesian linear regression to compare the impact of the different shrinkage prior distributions on the regression coefficients. Specifically, we compare the Gambel distribution (with various hyperparameter settings) to the Horseshoe, Laplace, and spike and slab priors. Table \ref{tab:priors} details the prior distributions and their hyperparameters.

For a fair comparison, the different prior distributions need to be calibrated. Typically, this involves matching the mean and variance. Since some priors lack moments, we need to use an alternative criterion and, instead, examine the curvature of the density function. We calibrate the priors by ensuring that they have a similar mass in a ``modal region" around the origin.  Recently, for unimodal distributions, \citet{cabral2024does} proposed to delimit the bulk of the distribution from
its tails, using the point of
maximum curvature. Intuitively, the curvature is the amount by
which a curve deviates from being a straight line. This implies that the
 modal region is defined as the interval $(-\text{PMCurv}, \text{PMCurv})$, where 
$$\text{PMCurv} := \argmax_{x>0} \frac{f''(x)}{(1+f'(x)^2)^{3/2}}$$
is the Point of Maximum Curvature,
which can be determined numerically. 
\begin{table}[]
\centering
\caption{Prior distributions, whether or not they have a singularity at 0, Point of Maximum Curvature (PMCurv) and $P(|\theta| < \text{PMCurv})$. }
\resizebox{\columnwidth}{!}{%
\begin{tabular}{lccc}
\toprule
Prior                                         & \multicolumn{1}{l}{Singularity?} & \multicolumn{1}{l}{PMCurv} & \multicolumn{1}{l}{$P(|\theta| < \text{PMCurv})$} \\ \midrule
Horseshoe $(q=2, a = 0.5, b = 0.5, \xi = 1)$ & Yes                              & 0.195                      & 0.27                                               \\
Gambel 1 $(q=2, a = 0.5, b = 0.52, \xi = 1)$   & No                               & 0.183                      & 0.25                                               \\
Gambel 2 $(q=2, a = 0.3, b = 1.6, \xi = 1)$    & No                               & 0.98                       & 0.27                                               \\
Laplace $(\lambda = 1)$                       & No                               & 0.34                       & 0.29                                               \\
Spike and slab                                & Yes                              & -                          & -   \\                               \bottomrule               
\end{tabular}}
\label{tab:priors}
\end{table}

As shown in Table \ref{tab:priors}, we chose the hyperparameters of the various priors so that their masses at the modal region, $P(|\theta|<\text{PMCurv})$, are approximately the same. The spike and slab prior is defined by
$\theta_i \overset{i.i.d.}{\sim} \pi N(0, \sigma_\theta^2) +  (1-\pi) \delta_ 0 $
where $\sigma_\theta^2 = 1000$ and $\pi \sim Beta(2.85, 1)$. This hyperparameter configuration ensures that the mass at zero is similar to the mass in the modal region of the other priors. Specifically, we have $E[\pi] = P(|G| < \text{PMCurv}_G) = 0.27$ and $SD[\pi] = 0.2$. 

\begin{figure}[]
    \centering
    \includegraphics[width=0.39\textwidth]{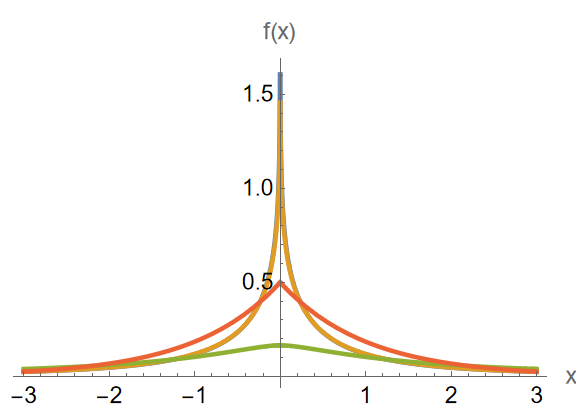} 
    \includegraphics[width=0.6\textwidth]{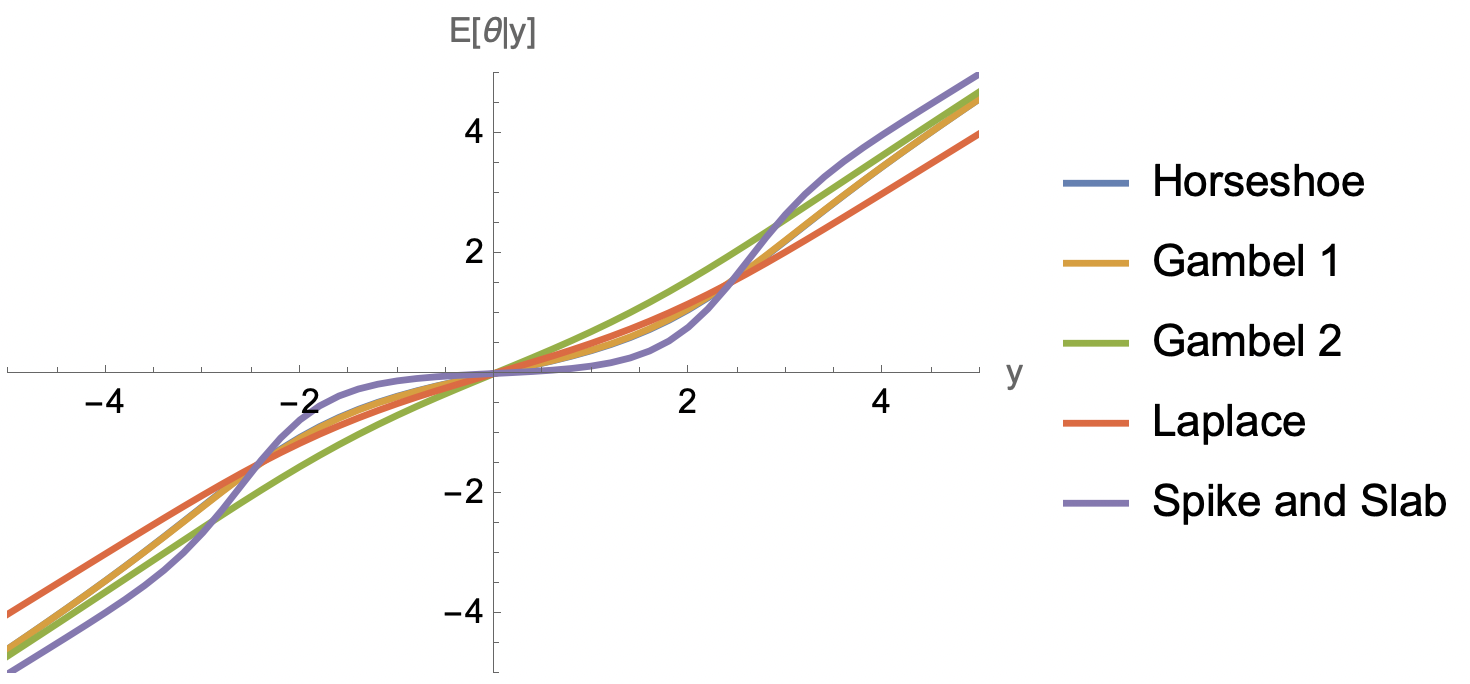} 
    \caption{Prior densities (left) and their shrinkage profiles (right). The shrinkage profile of the Horseshoe prior cannot be seen as it is superimposed by the Gambel 1 prior.}
    \label{fig:shrink}
\end{figure}

The density and shrinkage profiles of the priors are shown in Figure \ref{fig:shrink}. The shrinkage profiles contrast the posterior expectation of the regression coefficient, $E[\theta|y]$, with a varying datapoint $y$, where the likelihood is $y \sim N(\theta,1)$, and different priors are chosen for $\theta$. We note that Gambel 1 is very close to the Horseshoe prior, while Gambel 2 is closer to a uniform prior (see Table \ref{tab:priors} and Figure \ref{fig:shrink}). 

To evaluate the performance of the shrinkage priors across various scenarios, we conduct simulations considering different levels of sparsity in the regression coefficients, variability in the data, and sample sizes. The posterior inference is performed using three separate Markov Chains Monte Carlo (MCMC) chains implemented in the software JAGS \citep{plummer2003jags}, where each chain is run for 10,000 iterations, with the first 3,500 iterations discarded as burn-in.

\subsection{Simulation study 1}

In this simulation, we focus on the case $p<n$ (fewer regression coefficients than datapoints) and consider the linear regression model
$\mathbf{y}\mid\mathbf{X},\boldsymbol{\theta},\sigma \sim N(\mathbf{X}\boldsymbol{\theta}, \sigma^2\mathbf{I}_N)$,
with data simulated using the following parameter value configurations:
\begin{enumerate}
    \item $\boldsymbol{\theta}^{\text{sim}} = (20, 20, 40, 40, 0, 0, 0, 0)'$,  $\sigma=3$, $N=50$. The design matrix $\mathbf{X}$ is drawn from a normal distribution with mean 0, marginal variance 1, and the pairwise correlation between predictors is 0.5.  
    \item Same as before, but $\boldsymbol{\theta}^{\text{sim}} = (1, 1, 20, 20, 0, 0, 0, 0)'$.
    \item Same as before, but $\boldsymbol{\theta}^{\text{sim}} = (\underbrace{20,\ldots,20}_{\text{7 times}}, \underbrace{5,\ldots,5}_{\text{8 times}}, \underbrace{0,\ldots,0}_{\text{15 times}})'$, $\sigma=15$, $N=300$. 
    \item  $\boldsymbol{\theta}^{\text{sim}} = (\underbrace{20,\ldots,20}_{\text{7 times}}, \underbrace{5,\ldots,5}_{\text{8 times}}, \underbrace{0,\ldots,0}_{\text{15 times}})',$ and $\boldsymbol{x}_j = Z_1 + w_j$, $j=1,\ldots,5;$ $\boldsymbol{x}_j = Z_2 + w_j$, $j=6,\ldots,10;$ $\boldsymbol{x}_j = Z_3 + w_j$, $j=11,\ldots,15;$ and $\boldsymbol{x}_j \sim N(0,1)$, $j=16,\ldots,30$. Here, $Z_1,Z_2,Z_3$ are i.i.d. standard normal variables and $w_j\sim N(0,0.01)$. We have $N=300$.
    \item Same as before, but $\boldsymbol{\theta}^{\text{sim}} = (\underbrace{1,\ldots,1}_{\text{7 times}}, \underbrace{50,\ldots,50}_{\text{8 times}}, \underbrace{0,\ldots,0}_{\text{15 times}})'$, $\sigma=15$, $N=500$.
    \item  $\boldsymbol{\theta}^{\text{sim}} = (\underbrace{0.5,\ldots,0.5}_{\text{5 times}}, \underbrace{5,\ldots,5}_{\text{5 times}}, \underbrace{40,\ldots,40}_{\text{7 times}}, \underbrace{0,\ldots,0}_{\text{15 times}})'$, $\sigma=15$, $N=300$, and design matrix of same structure as before.
\end{enumerate}
These six configurations are adapted from \cite{tibshirani1996regression}, \cite{zou2005regularization}, \cite{kyung2010penalized}, \cite{li2010bayesian}, and \cite{roy2017selection}. For each, we generated 100 replicates. We evaluate  the posterior means $\hat{\boldsymbol{\theta}} = E[\boldsymbol{\theta}|\mathbf{y}]$. The code of the JAGS model considering the Gambel priors is shown in the Supplement. We show in Table \ref{tab:priors2} the median of the realised squared error loss $\tfrac{1}{N}\sum_{i=1}^N (\hat{\theta}_i - \theta_i^{\text{sim}})^2$.

\begin{table}[]
\centering
\caption{Median of  squared error loss for different prior distributions and configurations in the simulation study 1.}
\begin{tabular}{lcccccc}
\toprule
               & Conf. 1        & Conf. 2        & Conf. 3        & Conf. 4        & Conf. 5        & Conf. 6        \\ \midrule
Horseshoe      & 0.247          & 0.274          & 1.165          & \textbf{3.814} & \textbf{1.519} & \textbf{2.146} \\
Gambel 1       & 0.246          & 0.272          & 1.166          & 3.887          & 1.528          & 2.170          \\
Gambel 2       & \textbf{0.187} & \textbf{0.238} & 1.483          & 7.831          & 1.839          & 3.833          \\
Laplace        & 0.328          & 0.323          & 1.270          & 14.023         & 4.955          & 8.844          \\
Spike and slab & 0.172          & 0.320          & \textbf{1.125} & 3.912          & 1.696          & 2.492 \\
\bottomrule
\end{tabular}
\label{tab:priors2}
\end{table}

Unsurprisingly, the Horseshoe and Gambel 1 performed similarly. Thus, the fact that the Gambel 1 prior, unlike the Horseshoe, has no singularity at zero, little affects the results. The simulations also show that the Horseshoe and Gambel 1 prior performed similarly or better than the Spike and slab prior, and overall, the Laplace prior performed worst. The superiority of the Horseshoe over the Laplace prior has previously been noted by  \cite{carvalho2010horseshoe}.

\subsection{Simulation study 2}

In this simulation, we consider data $\mathbf{y}^{\text{sim}} = (y_1^{\text{sim}},\dotsc, y_N^{\text{sim}})$ generated as
$y_i^{\text{sim}} \sim N(\theta_i^{\text{sim}},1)$, with $\theta_i^{\text{sim}} \sim w t_\nu(0,3) + (1-w)\delta_0$ for each $i$, 
where $t_\nu(0,3)$ is a centred Student-t distribution with $\nu$ degrees of freedom and scale 3.  We set $N=1000$, $\epsilon \in \{2, 10\}$ and $w \in \{0.05, 0.2, 0.5\}$. Smaller $w$ leads to sparser simulated signals. The smaller $\nu$, the larger the kurtosis of the non-zero regression coefficients (more values closer to zero with potentially more outliers). 
A similar simulation scenario is considered by \cite{carvalho2010horseshoe}.
We fit a simple linear regression model $\mathbf{y} \sim N(\boldsymbol{\theta},\tau^{-1}\mathbf{I}_N)$, where $\tau$ has a truncated Cauchy distribution and each $\theta_i$ follows independently one of the priors in Table \ref{tab:priors}. Table \ref{tab:priors1} shows the realised squared error losses.

\begin{table}[]
\centering
\caption{Squared error loss for different prior distributions in the simulation study 2.}
\begin{tabular}{lcccccc}
                \toprule
               & \multicolumn{3}{c}{$\epsilon = 2$} & \multicolumn{3}{c}{$\epsilon = 10$} \\
               \cmidrule(lr){2-4}\cmidrule(lr){5-7}
               & $w=0.05$   & $w=0.2$   & $w=0.5$   & $w=0.05$    & $w=0.2$   & $w=0.5$   \\ \midrule
Horseshoe      & 0.519                       & 0.535                       & 0.572                       & 0.530                       & 0.535                       & \textbf{0.589}              \\
Gambel 1       & 0.510                       & 0.529                       & \textbf{0.571}                       & 0.522                       & 0.527                       & 0.590                       \\
Gambel 2       & \textbf{0.211}              & 0.334              & 1.513                       & 0.197                       & 0.350              & 1.189                       \\
Laplace        & 1.133                       & 0.670                       & 1.632                       & 0.734                       & 0.678                       & 0.919                       \\
Spike and slab & 0.408                       & \textbf{0.282}                       & 0.950              & \textbf{0.167}              & \textbf{0.324}                       & 1.124  \\
\bottomrule
\end{tabular}
\label{tab:priors1}
\end{table}

Again, the Horseshoe and Gambel 1 priors led to very similar performances. They also gave better results than the other priors when the signal is sparser ($w=0.5$). Overall, the Laplace prior performed worst, and the Gambel 2 prior (which is closest to a uniform prior) tended to give better results for non-sparse signals.    

\section{Conclusion}\label{sec:discussion}

The results of this paper provide valuable insights into the properties and applications of the newly proposed Gambel distribution, a generalisation of many shrinkage priors. By placing the analysis in reliability theory and wealth distributions, we extend understanding beyond traditional Horseshoe-type priors. Empirical evaluations show that the Gambel distribution, particularly its first variant (Gambel 1), performs comparably to the Horseshoe prior in shrinkage efficiency under high sparsity without a singularity at 0. Its robustness across different sparsity levels and signal-to-noise ratios highlights its practical utility.

We also examined tail behaviour and hazard rates, linking the Gambel distribution to wealth distributions in order to quantify inequality and concentration effects via Lorenz curve and Gini index analyses. This dual perspective bridges gaps in the literature by addressing behaviour at both the origin and the tails, offering applications in high-dimensional regression with mixed sparse and dense signals.

Future research will explore the use of the Gambel prior in other regression models, develop efficient computational methods such as improved MCMC schemes, and assess performance on real-world datasets. These efforts will further establish the Gambel distribution's potential and encourage its adoption in statistics and machine learning.
%\bibliographystyle{chicago}
%\bibliography{prior_references}

%%%%%%%%%%%%%%%%%%%%

\newpage
\setcounter{equation}{0}
\setcounter{table}{0}
\setcounter{figure}{0}
\setcounter{section}{0}
\numberwithin{table}{section}
\renewcommand{\theequation}{S.\arabic{equation}}
\renewcommand{\thesection}{S.\arabic{section}}
\renewcommand{\thesubsection}{S.\arabic{section}.\arabic{subsection}}
\renewcommand{\thetable}{S.\arabic{table}}
\renewcommand{\thefigure}{S.\arabic{figure}}
\begin{center}
\textbf{\Large{Supplement }}\\

\end{center}
 \section{Proofs}
 \subsection{Proof of Theorem~1} \label{proof:theo1}
 \begin{proof}
The conditional distribution of $\theta$ given $\kappa$ is an Exponential Power distribution with parameters $q$ and $\kappa$, i.e. it has density 
\[
\text{EP}(\theta \mid q,\kappa) = \frac{q\,\phi(q)}{(1/\kappa - 1)^{1/q} 2\Gamma(1/q)}\exp\left\{\frac{-\kappa (\phi(q)\, |\theta|)^q}{1-\kappa}\right\} \qquad \text{where }\phi(q) = \left(\sftwo{\Gamma}{3/q}{1/q}{}\right)^{1/2} 
\]
and the marginal distribution of $\kappa$ is a Generalized Three Parameters Beta, i.e. it has a density
\[
\text{G3B}(\kappa \mid a, b, \xi) = \frac{1}{B(a,b)} \frac{\xi^a \kappa^{a-1} (1-\kappa)^{b-1}}{(1-(1-\xi)\kappa)^{a+b}}
\]
Thus, setting 
\[
C := \frac{q\,\phi(q)\, \xi^a}{ 2\Gamma(1/q)\, B(a,b)}
\]
the marginal density for $\theta$ is 
\begin{align*}
f(\theta\mid q, a, b, \xi) =& \,C\, \int_0^1 \frac{1}{(1/\kappa - 1)^{1/q}} \exp\left\{\frac{-\kappa (\phi(q)\, |\theta|)^q}{1-\kappa}\right\} \frac{\kappa^{a-1} (1-\kappa)^{b-1}}{(1-(1-\xi)\kappa)^{a+b}} \dd \kappa\\
=& \,C\, \int_0^1 \exp\left\{-\alpha_q(\theta)\frac{\kappa}{1-\kappa}\right\}\kappa^{a + 1/q - 1} (1-\kappa)^{b - 1/q - 1} (1 - (1-\xi)\kappa)^{-(a+b)} \dd \kappa \\
&\text{where } \alpha_q(\theta) := (\phi(q)\, |\theta|)^q
\end{align*}
Substituting $y:= \frac{\kappa}{1-\kappa}$, we have $\kappa = \frac{y}{1+y}$, $\text{d} \kappa = \frac{\text{d} y}{(1+y)^2} $, and 
\begin{align*}
f(\theta\mid q, a, b, \xi) =& \,C\, \int_0^1\,\exp\left\{- \alpha_q(\theta)\,y\right\}\left(\frac{y}{1+y}\right)^{a + 1/q - 1} (1+y)^{-b + 1/q + 1} \times \\
&\times \left(1 - (1-\xi)\frac{y}{1+y}\right)^{-(a+b)} (1+y)^{-2}\dd y\\
=&\, C\, \int_0^1\,\exp\left\{-\alpha_q(\theta)\,y\right\} y ^{a+1/q-1} (1 + \xi\,y)^{-(a+b)} \dd y
\end{align*}
The latter integral can be found in \cite{gradshteyn2014table}, leading to%, 3.383(5),
\begin{align*}
f(\theta\mid q, a, b, \xi) =& \,C\, \xi^{-a - 1/q}\Gamma(a + 1/q) \Psi\left( a+\tfrac{1}{q},1+\tfrac{1}{q}-b, \tfrac{(\phi(q)\,|\theta|)^q}{\xi}\right)
\end{align*}
\end{proof}

\subsection{Proof of Corollary~1}
\begin{proof}
Since
\[
    \Psi\left( x,y,z\right) = \frac{\Gamma(1-y)}{\Gamma(1+x-y)}M(x,y,z) + \frac{\Gamma(y-1)}{\Gamma(x)}z^{1-y}M(1+x-y,2-y,z), \text{for } y\notin \mathbb{Z}
\]
it follows that $\lim_{\theta\rightarrow 0} \Gambel(\theta\mid q,a,b,\xi) = \lim_{\theta\rightarrow 0}\frac{\Gamma(1/q-b)}{\Gamma(a+1/q)}z^{b-1/q}$, from which it follows the behaviour at the origin. 
Contrary, using the series expansion at $z = \infty$ of the Tricomi function, we get that $ \Gambel(\theta\mid q,a,b,\xi) \sim (\phi(q) |\theta|)^{-aq - 1} / (\xi^{a+1/q})$.
\end{proof}

 \subsection{Proof of Proposition~1}
\begin{proof}
Note that if $\kappa\sim\text{G3B}(\kappa\mid a,b,\xi)$, then there exists $X_0$ and $X_1$ independent random variables such that $X_1 / (X_0 + X_1) = \kappa$. In particular, 
\[
X_0 \sim \text{Gamma}(X_0 \mid b,\beta_0) \qquad X_1 \sim \text{Gamma}(X_1 \mid a,\beta_1)
\]
where the Gamma distribution is parameterized in terms of the shape and rate parameters, meaning that $\mathbb{E}[X_0] = \frac{b}{\beta_0}$. Moreover, the rate parameters $\beta_0$ and $\beta_1$ are such that $\xi = \frac{\beta_1}{\beta_0}$. The proof of this equivalence follows directly from the construction of the \text{G3B} in \cite{libby1982multivariate}.

 Thus, it follows that $\lambda^2 =\frac{1-\kappa}{\kappa} = \frac{X_0}{X_1}$ and
 \[
 \lambda^2 \sim \text{GB2}(\lambda^2 \mid 1, \xi, b, a) 
 \]
 since the ratio of two independent Gamma random variables is distributed according to a generalized beta prime distribution.
\end{proof}

 \subsection{Proof of Theorem~2}
\begin{proof}
The density of $X$ and $Y$ are respectively 
\[
\text{GG}(x \mid \xi^{1/q}, 1, q) = \frac{q}{\xi^{1/q} \Gamma(1/q)}\exp\left\{-\left(\frac{x}{\xi^{1/q}}\right)^q\right\}
\]
and 
\[
\text{GB2}(y\mid q, \phi(q), a, b) = \frac{q (y/\phi(q))^{aq-1} (1+(y/\phi(q))^q)^{-(a+b)}}{\phi(q) B(a,b)}
\]
Thus the density of $Z = X/Y$ is 
\begin{align*}
    f(z\mid q, a, b, \xi) = \int_{0}^{+\infty} \text{GG}(zy \mid \xi^{1/q}, 1, q) \text{GB2}(y\mid q, \phi(q), a, b) |y| \dd y 
\end{align*}
Setting 
\begin{align*}
C := \frac{q^2}{\xi^{1/q} \Gamma(1/q) B(a,b)}
\end{align*}
the density of $Z$ is 
\begin{align*}
    f(z\mid q, a, b, \xi) = C\,\int_{0}^{+\infty} 
    \exp\left\{-\left(\frac{zy}{\xi^{1/q}}\right)^q\right\}
   \left(\frac{y}{\phi(q)}\right)^{aq} \left(1+\left(\frac{y}{\phi(q)}\right)^q\right)^{-(a+b)} \dd y 
\end{align*}
Substituting $x:= \left(\frac{y}{\phi(q)}\right)^q$, we have $y = \phi(q) x^{1/q}$, $\dd y = \frac{\phi(q)}{q}x^{1/q-1} \dd x$, so 
\begin{align*}
  f(z\mid q, a, b, \xi) = C\,\frac{\phi(q)}{q}\int_{0}^{+\infty} 
    \exp\left\{-\left(\frac{z\phi(q)}{\xi^{1/q}}\right)^q x\right\}
   x^{a+1/q-1} \left(1+x\right)^{-(a+b)} \dd x
\end{align*}
the integral above can be found in \cite{gradshteyn2014table}, 3.383(5),
\begin{align*}
 f(z\mid q, a, b, \xi) = C\,\frac{\phi(q)}{q}\Gamma\left(a +\frac{1}{q}\right)\Psi\left(a+\frac{1}{q}, 1 +\frac{1}{q} - b, \frac{(\phi(q) z)^q}{\xi}\right)
\end{align*}
\end{proof}

 \subsection{Proof of Proposition~2}
\begin{proof}
Employing the ratio representation provided by Theorem~\ref{thm:ratio}, i.e. $\theta = X/Y$ such that $X\sim \text{GG}(\xi^{1/q},d=1,q)$ and $Y\sim GB2(q,\phi(q),a,b)$, with $X \perp Y$. 
The c.d.f. for $X$ is 
\[
F_{\text{GG}}(x \mid \xi^{1/q}, 1, q) = \frac{\gamma(1/q, x^q / \xi)}{\Gamma(1/q)}
\]
where $\gamma$ is the lower incomplete gamma function: $\gamma(1/q,  x^q / \xi) = \int_{0}^{ x^q / \xi}t^{1/q}e^{-t}\dd t$.
\begin{align*}
F(\theta; q,a,b,\xi) &= \mathbb{P}(X/Y \leq \theta) = \mathbb{E}[\mathbb{P}(X \leq \theta\,Y\mid Y)]=
\mathbb{E}\left[\frac{\gamma(1/q, (\theta Y)^q/\xi)}{\Gamma(1/q)}\right]\\
&= \int_{0}^{+\infty} \frac{\gamma(1/q, (\theta y)^q/\xi)}{\Gamma(1/q)} \frac{q (y/\phi(q))^{aq-1} (1+(y/\phi(q))^q)^{-(a+b)}}{\phi(q) B(a,b)}\dd y
\end{align*}
Defining 
\[
C:= \frac{q}{\Gamma(1/q)\phi(q) B(a,b)}
\]
we have 
\begin{align*}
F(\theta; q,a,b,\xi) &= C \int_{0}^{+\infty} \gamma\left(\frac{1}{q}, \frac{(\theta y)^q}{\xi}\right) 
\left(\frac{y}{\phi(q)}\right)^{aq-1} \left(1+\left(\frac{y}{\phi(q)}\right)^q\right)^{-(a+b)}\dd y
\end{align*}
Substituting $x:= y^q$, we have $y = x^{1/q}$, $\dd y = \frac{x^{1/q-1}}{q} \dd x$, so 
\begin{align*}
F(\theta; q,a,b,\xi) &=  \frac{\phi(q)^{qb}}{\Gamma(1/q)B(a,b)} \int_{0}^{+\infty} 
\gamma\left(\frac{1}{q}, \frac{\theta^q}{\xi} x\right) 
x^{a-1}
\left(\phi(q)^{q}+\, x\right)^{-(a+b)}\dd x
\end{align*}
the integral above can be found in \cite{prudnikov1988integrals}, 2.10.2(4),
\begin{multline*}
F(\theta; q,a,b,\xi) =\left(\tfrac{\phi(q)\,\theta}{\xi^{1/q}}\right)
\tfrac{  q\,\Beta\left(a+\tfrac{1}{q}, b-\tfrac{1}{q}\right) }{\Gamma\left(\tfrac{1}{q}\right)\Beta(a,b)}
\sffour{\, _2F_2}{\tfrac{1}{q},}{ a +\tfrac{1}{q} ;}{ 1+\tfrac{1}{q},}{ 1+\tfrac{1}{q} -a;}{ \left(\tfrac{\phi(q)\,\theta}{\xi^{1/q}}\right)^q  }  \\
+\left(\tfrac{\phi(q)\,\theta}{\xi^{1/q}}\right)^{qb}\tfrac{\Gamma\left(\tfrac{1}{q} - b\right) }{b \Gamma\left(\tfrac{1}{q}\right)\Beta({b},{a})} 
\sffour{\,  _2F_2}{b,}{ b+a;}{ 1+b,}{ 1+b-\tfrac{1}{q};}{\left(\tfrac{\phi(q)\,\theta}{\xi^{1/q}}\right)^q }
\end{multline*}
\end{proof}

\subsection{Proof of Proposition~3}
\begin{proof}
Employing the ratio representation provided by Theorem~\ref{thm:ratio}, i.e. $\theta = X/Y$ such that $X\sim \text{GG}(\xi^{1/q},d=1,q)$ and $Y\sim GB2(q,\phi(q),a,b)$, with $X \perp Y$, note that by the independence, $\mathbb{E}	((X/Y)^k) =  \mathbb{E}(X^k) \mathbb{E}(1/Y^{k})$ with
\begin{align*}
\mathbb{E}(X^k)				& =  \xi^{k/q}\sftwo{\Gamma}{(d+k) / q}{d/q}{} =  \xi^{k/q}\sftwo{\Gamma}{(1+k) / q}{1/q}{}\\ 
\mathbb{E} 	(1/Y^k) 		& =  \int_0^{\infty}  \frac{q}{\phi(q)^{aq}\Beta(a,b)} 
\frac{y^{aq-k-1}}{(1+(y/\phi(q))^q)^{b+a}}   \id y \\
&=  \int_0^{\infty} \frac{q}{\phi(q)^{aq}\Beta(a,b)}
\frac{\phi(q)}{q}\phi(q)^{aq -k-1}w^{a-k/q-1}\frac{1}{(1+w)^{b+a}}
  \id w  \quad[\text{Subst } w=(y/\phi(q))^q] \\
& = \frac{1}{\phi(q)^k\Beta(a,b)} 
\underbrace{\int_{0}^{\infty}\frac{w^{a-k/q-1}}{(1+w)^{a-k/q+k/q+b}}\id w}_{\text{kernel of Beta prime distribution}}\\
& = \frac{\Beta(a-k/q,b+k/q)}{\phi^k\Beta(a,b)}
\end{align*}

For the approximation of the kurtosis use $\Gamma(x+a)/\Gamma(x+b)\sim x^{a-b}$, and note
\begin{align*}
\text{Kurtosis}&\sim \left(\tfrac{1}{q}\right)^{-2/q} 
\left(\tfrac{3}{q}\right)^{2/q} 
(b-\tfrac{4}{q})^{-2/q} a^{2/q} (b+\tfrac{2}{q})^{2/q} b^{-2/q} \\ 
& = 3^{2/q}  \left(\frac{a(qb+2)}{b(qa-4)}\right)^{2/q}
\end{align*}

\end{proof}

 \subsection{Proof of Theorem~3}
\begin{proof}
Note that for any hazard function, the following expression for the first derivative holds true:
\[
r'(x) = \frac{1}{S(x)}\left(\frac{[S'(x)]^2}{S(x)} - S''(x)\right)
\]
where, for all $x>0$, $S(x) \geq 0$ and $S'(x) \leq 0 $. 
Thus, the hazard function $r(x)$ is decreasing for those $x$ such that
\[
\frac{S''(x)}{S'(x)}<\frac{S'(x)}{S(x)}
\]
Considering the survival function, we have 
\begin{align*}
S(x) = 1 -\left(\tfrac{\phi(q)\,x}{\xi^{1/q}}\right)
\tfrac{  q\,\Beta\left(a+\tfrac{1}{q}, b-\tfrac{1}{q}\right) }{\Gamma\left(\tfrac{1}{q}\right)\Beta(a,b)}
\sffour{\, _2F_2}{\tfrac{1}{q},}{ a +\tfrac{1}{q} ;}{ 1+\tfrac{1}{q},}{ 1+\tfrac{1}{q} -a;}{ \left(\tfrac{\phi(q)\,x}{\xi^{1/q}}\right)^q  }  \\
-\left(\tfrac{\phi(q)\,x}{\xi^{1/q}}\right)^{qb}\tfrac{\Gamma\left(\tfrac{1}{q} - b\right) }{b \Gamma\left(\tfrac{1}{q}\right)\Beta({b},{a})} 
\sffour{\,  _2F_2}{b,}{ b+a;}{ 1+b,}{ 1+b-\tfrac{1}{q};}{\left(\tfrac{\phi(q)\,x}{\xi^{1/q}}\right)^q }
\end{align*}
where, for $x\rightarrow\infty$, \citep[see,][]{paris2005kummer}
\[
\sffour{\,  _2F_2}{a,}{ d;}{ b,}{c;}{x} \sim 
\frac{\Gamma(b) \Gamma(c)}{\Gamma(a) \Gamma(d)} x^{(a+d) - (b+c)} e^{x}
\]
Thus,
\begin{align*}
S(x) \sim 1 - e^{\left(\tfrac{\phi(q)\,x}{\xi^{1/q}}\right)^q } \bigg\{  x^{2qa - 2q + 1} \, \left(\tfrac{\phi(q)}{\xi^{1/q}}\right)^{2qa - 2q + 1}
\tfrac{  q\,\Beta\left(a+\tfrac{1}{q}, b-\tfrac{1}{q}\right) }{\Gamma\left(\tfrac{1}{q}\right)\Beta(a,b)}
\frac{\Gamma( 1+\tfrac{1}{q}) \Gamma(1+\tfrac{1}{q} -a)}{\Gamma(\tfrac{1}{q}) \Gamma(a +\tfrac{1}{q})}
 \\
+ x^{q(a+b) - 2q + 1} \left(\tfrac{\phi(q)}{\xi^{1/q}}\right)^{q(a+b) - 2q + 1}\tfrac{\Gamma\left(\tfrac{1}{q} - b\right) }{b \Gamma\left(\tfrac{1}{q}\right)\Beta({b},{a})} 
\frac{\Gamma(1+b) \Gamma(1+b-\tfrac{1}{q})}{\Gamma(b) \Gamma(b+a)}  \bigg\}.
\end{align*}

Moreover, using the series expansion at $z=\infty$ of the Tricomi function have that 
\[
\Psi(a,b,x) \sim x^{-a} 
\]
and thus
\[
S'(x) \sim - x^{-aq-1}\,\frac{\Gamma(a+ 1/q)}{\Gamma(1/q)} \frac{q}{B(a,b)} \left(\frac{\phi(q)}{\xi^{1/q}}\right)^{-aq}
\]

Finally, note that \citep[see, for instance][]{abramowitz1948handbook}
\[
\Psi(a,b,x) = x^{-a}  {}_2F_0(a, 1 + a - b;; -x^{-1})
\]
so 
\[
\frac{S'(x)}{S''(x)} = -  \frac{x}{aq + 1} 
\frac{ {}_2F_0\left(a + 1/q , a - b;; -\left(\frac{\phi(q)\,x}{\xi^{1/q}}\right)^{-q}\right)}
{ {}_2F_0\left(a + 1/q + 1, a - b;; -\left(\frac{\phi(q)\,x}{\xi^{1/q}}\right)^{-q}\right)} \sim -  \frac{x}{aq + 1}
\]

\end{proof}

 \subsection{Proof of Theorem~4}
\begin{proof}

We begin with the expression
$$
L(p) = \frac{\int\limits_{0}^{F^{-1}(p)}x f(x) \text{d} x}{\int\limits_{0}^{+\infty}x f(x) \text{d} x}
$$
This can be rewritten as
$$L(p) = \frac{1}{\mu} \int\limits_0^{F^{-1}(p)} \theta \, \text{Gambel} (\theta \mid q,a,b,\xi) \text{d} \theta$$

Using the ratio representation from Theorem~\ref{thm:ratio}, we obtain
$$L(p) = \frac{C\phi(q)}{\mu\,q} \int\limits_0^{F^{-1}(p)} \int\limits_{0}^{+\infty} \theta \exp\left\{-\frac{\phi(q)^q}{\xi} \theta^q x\right\} x^{a+1/q-1} (1+x)^{-(a+b)}\text{d}x\text{d}\theta$$
where
$$C = \frac{10 q^2 \xi ^{-1/q}}{\Gamma \left(\frac{1}{q}\right)}, \quad \mu = \frac{\xi ^{1/q} \, _3F_3\left(1-\frac{1}{q},10+\frac{1}{q},\frac{2}{q};1,10,\frac{1}{q};1\right)}{\phi(q)}$$

Applying Fubini's theorem, we can interchange the order of integration
$$L(p) = \frac{C\phi(q)}{\mu\,q} \int\limits_{0}^{+\infty} \left( x^{a+1/q-1} (1+x)^{-(a+b)}\int\limits_0^{F^{-1}(p)} \theta \exp\left\{-\frac{\phi(q)^q}{\xi} \theta^q x\right\}\text{d}\theta\right)\text{d}x$$

Now, let us evaluate the inner integral
$$\int\limits_0^{F^{-1}(p)} \theta \exp\left\{-\frac{\phi(q)^q}{\xi} \theta^q x\right\}\text{d}\theta = \frac{\left(\frac{x \phi(q)^q}{\xi }\right)^{-2/q} \left(\Gamma \left(\frac{2}{q}\right)-\Gamma \left(\frac{2}{q},\frac{x \left(\phi(q) \theta ^*\right)^q}{\xi }\right)\right)}{q}$$
where $\Gamma (a,z)=\int _z^{\infty }d t t^{a-1} e^{-t}$ is the incomplete Gamma function.

Substituting this result back into the outer integral and evaluating, we obtain

\begin{align*}
&\int\limits_{0}^{+\infty} \left( x^{a+1/q-1} (1+x)^{-(a+b)}\int\limits_0^{F^{-1}(p)} \theta \exp\left\{-\frac{\phi(q)^q}{\xi} \theta^q x\right\}\text{d}\theta\right)\text{d}x \\ 
&=\frac{\left(\theta ^*\right)^{b q+1} \xi ^{\frac{1}{q}-b} \Gamma \left(\frac{1}{q}-b\right) \phi(q)^{b q-1} \, _2F_2\left(a+b,b+\frac{1}{q};b-\frac{1}{q}+1,b+\frac{1}{q}+1;\frac{\phi(q)^{\frac{q}{2}} \left(\theta ^*\right)^q}{\xi }\right)}{b q+1} \\ &+\frac{\left(\theta ^*\right)^2 \Gamma \left(a+\frac{1}{q}\right) \Gamma \left(b-\frac{1}{q}\right) \, _2F_2\left(a+\frac{1}{q},\frac{2}{q};-b+\frac{1}{q}+1,1+\frac{2}{q};\frac{\phi(q)^{\frac{q}{2}} \left(\theta ^*\right)^q}{\xi }\right)}{2 \Gamma (a+b)}
\end{align*}

Finally, multiplying this result by $(C\phi(q))/(\mu q)$ and simplifying yields the theorem.

\end{proof}

\subsection{Proof of Proposition~5}
\begin{proof}
A direct calculation shows
$\int\limits_{a}^{\infty}be^{-bu}\,\text{d}u = e^{-ba}$ for $a,b>0$.
Therefore, with $\alpha:=1+1/q$ and since $\Gamma(\alpha)=\Gamma(\frac{1}{q})/q$, %\todoB{Matching with the notation of the EP, $\gamma = ((1-\kappa) / \kappa)^{1/q}$ and $\eta = \phi(q)$}
\begin{align*}
\pi(\theta\mid\gamma)&=\frac{q\phi(q)}{2\Gamma(1/q)\gamma}
\exp\left\{-\left(\frac{\phi(q)|\theta|}{\gamma}\right)^q\right\}\\
&=\frac{\phi(q)}{2\Gamma(1+1/q)\gamma}
\int\limits_{u>\frac{|\theta|^q}{\gamma^q}}
\phi(q)^q e^{-\phi(q)^qu}\textrm{d}u \\
&=\frac{1}{\Gamma(1+1/q)}
\int\limits_{\gamma u^{1/q}>|\theta|}\frac{1}{2\gamma u^{1/q}}
\phi(q)^{q(1/q+1)} \exp\left\{-\phi(q)^qu\right\}u^{1/q+1-1}\textrm{d}u \\
&=\frac{1}{\Gamma(\alpha)}
\int\limits_{\gamma u^{1/q}>|\theta|}\frac{1}{2\gamma u^{1/q}}
\phi(q)^{q\alpha} \exp\left\{-\phi(q)^q u\right\}u^{\alpha-1}\textrm{d}u 
\end{align*}
\end{proof}

 \subsection{Proof of Theorem~5} \label{proof:consistency}

\begin{proof}
Recall that 
\[
\pi(\theta) 
= \tfrac{c_1}{\xi^{1/q} }
\Psi\left(a+\tfrac{1}{q}, \tfrac{1}{q}+1-b,\tfrac{(\phi(q) |\theta|)^q}{\xi}\right)
\quad\text{with } c_1= \tfrac{\phi(q) q \Gamma(a+1/q)}{2\Gamma(1/q)\Beta(a,b)}
\]
For sufficiently large $n$, therefore

\begin{align*}
\Prob&\left(\|\bbeta_n-\bbeta_n^*\|
<\tfrac{\Delta}{n^{\rho/2}}\right) \\
&\overset{(a)}{\geq}
\left(1-\tfrac{p_n n^{\rho}\mathbb{E}(\theta_{nj}^2)}{\Delta^2}\right)  \left(\tfrac{\Delta}{(p_n n^{\rho})^{1/2}}\right)^{s_n} 
\left\{
\tfrac{c_1}{\xi^{1/q} }
\Psi\left(a+\tfrac{1}{q}, \tfrac{1}{q}+1-b,\tfrac{\phi(q)^q}{\xi}{\left[\|\bbeta^*_n\|_{\infty} +\tfrac{\Delta}{(p_n n^\rho)^{1/2}}\right]^q}\right)
\right\}^{s_n}\\
&\overset{(b)}{\geq}  \left(1-\tfrac{p_n n^{\rho}\mathbb{E}(\theta_{nj}^2)}{\Delta^2}\right)  \left(\tfrac{\Delta}{(p_n n^{\rho})^{1/2}}\right)^{s_n} 
\left\{
\tfrac{c_1}{\xi^{1/q} }
\Psi\left(a+\tfrac{1}{q}, \tfrac{1}{q}+1-b,\tfrac{\phi(q)^qc_q}{\xi}\left[\|\bbeta^*_n\|_{\infty}^q +\tfrac{\Delta^q}{(p_n n^\rho)^{q/2}}\right]\right)
\right\}^{s_n}\\
&\overset{(c)}{\gtrsim}
\left(1-\tfrac{p_n n^{\rho}\mathbb{E}(\theta_{nj}^2)}{\Delta^2}\right)  \left(\tfrac{\Delta}{(p_n n^{\rho})^{1/2}}\right)^{s_n} 
\left\{
\tfrac{c_1}{\phi(q)^{aq+1}c_q^{a+1/q}\xi^{1/q} }
\left[\tfrac{1}{\xi}\|\bbeta^*_n\|_{\infty}^q+\tfrac{\Delta^q}{(p_nn^{\rho})^{q/2}\xi}\right]^{-(a+1/q)}
\right\}^{s_n}
\end{align*}
where (a) follows from equation \eqref{eq:consistency_eq}.
For (b), we use that on $\R^2$ we have the relationship $\|\cdot\|_1\geq c_q\|\cdot\|_q$ between the $1$-norm and the $q$-(quasi)norm, where $c_q:=1$ for $q\geq 1$ and $c_q:=2^{q/(q-1)}$ for $q<1$ (the latter following from the generalized H\"{o}lder inequality $\|(u_iv_i)_{i=1}^N\|_{\tau}\leq\|u\|_p\|v\|_q$ for all $u,v\in\R^N$ whenever $\frac{1}{\tau}=\frac{1}{p}+\frac{1}{q}$ for $0<\tau\leq 1$ and $0<p,q<\infty$). Step (c) follows since  
$\Psi(\alpha, \beta, z) = z^{-\alpha}[1 + O(|z|^{-1})] \ \text{as} \ z \rightarrow +\infty$
\citep[see][13.1.8]{abramowitz1948handbook}.
(Indeed, for sufficiently large $z$, this implies 
$-\frac{M}{|z|}\leq\Psi(\alpha, \beta, z)z^{\alpha}-1\leq \frac{M}{|z|}$ for some $M>0$, and thus
$\Psi(\alpha, \beta, z)z^{\alpha}\geq 1 - \frac{M}{|z|} \geq  k$ for any $k\in(0,1)$ given that $z$ is sufficiently large; thus $\Psi(\alpha, \beta, z) \gtrsim z^{-\alpha}$.)

Taking the negative logarithm, it follows that for sufficiently large $n$:
\begin{align*}
-\log\Prob&\left(\|\bbeta_n-\bbeta_n^*\|
<\tfrac{\Delta}{n^{\rho/2}}\right) \\
&
\overset{(d)}{\lesssim} -\log\left( 1-\tfrac{p_n n^{\rho}c_2\xi^{2/q}}{\Delta^2}   
\right) 
 - s_n\log\left\{  \tfrac{\Delta  c_1}{(p_n n^{\rho})^{1/2}\xi^{1/q}\phi(q)^{aq+1}c_p^{a+1/q}}\right\}  \\
 &\hspace{6cm}
+ s_n\left(a+\tfrac{1}{q}\right) \log \left[\tfrac{1}{\xi} \|\bbeta^*_n\|_{\infty}^q+\tfrac{\Delta^q}{(p_n n^{\rho})^{q/2}\xi}\right] \\
&\overset{(e)}{=} -\log\left( 1-\tfrac{c_2C}{\Delta^2\log n}  
\right)  - s_n\log\left\{  \tfrac{\Delta c_1}{C^{1/2}\phi(q)^{aq+1}c_p^{a+1/q}}\right\}  
- \tfrac{s_n}{2} \log\log n   
\\ & \hspace{3cm}
+ \underbrace{s_n\left(a+\tfrac{1}{q}\right) \log \bigg[(\tfrac{1}{C}p_n n^{\rho}\log n)^{q/2}\|\bbeta^*_n\|_{\infty}^q}_{\textup{(}\star\textup{)}}+\Delta^q (\tfrac{1}{C}\log n)^{q/2}\bigg]
\end{align*}
where (d) follows from
Proposition \ref{prop:moments}, which guarantees under the assumption $\frac{2}{q}<a$ that $\mathbb{E}(\theta_{nj}^2)=c_2\xi^{2/q}$ for some constant $c_2$; and (e) follows since $\xi^{2/q} = C / (pn^{\rho} \log n)$. 
For $n \to \infty$, the dominating term in the final expression is ($\star$), and thus $-\log\Prob\left(\|\bbeta_n-\bbeta_n^*\| <\tfrac{\Delta}{n^{\rho/2}}\right) < dn$ for all $d>0$, which is the desired sufficient condition \eqref{StronConsistencySuffCond}. In more detail, ($\star$) can be rewritten as
$$
\left(a+\tfrac{1}{q}\right)\tfrac{q}{2}s_n\left\{\log p_n + \rho\log n + \log\log n + \tfrac{2}{q}\log \|\bbeta^*_n\|_{\infty}^q \right\}
$$
in which, by (L4) and (H4), $\|\bbeta^*_n\|_{\infty}^q $ is bounded as $n\to\infty$.
In the low-dimensional setting, the dominating term is thus proportional to $s_n\log n$, and the final result follows from condition (L3). %, $s_n\log n \prec n$. 
For the high-dimensional setting, the dominating term is proportional to $s_n\log p_n$ and (H3) guarantees the final result. 
\end{proof}

%%%%%%%%%%%%%%%%%%%%%%%%%%%%%%%%%%%%%%%%%%%%%%%%%%%%%%%%

\section{MCMC algorithms} \label{app:MCMC}

We next show the full conditionals of the random variables, which are necessary to implement the Gibbs sampler, first for the case $p < n$ and then for the case $p \geq n$.

\subsection{MCMC for the case $p < n$}

The full conditionals of the random variables in the hierarchical model in Section \ref{sec:MCMC1} are given next.

\begin{description}
\item[Update $\gamma_j$:] For $j=1,\ldots,p$
\begin{align*}
\pi(\gamma_j\mid \text{rest}) & \propto\gamma_j^{qa-1}
\exp\{-\lambda_j\gamma_j^q\} \tfrac{1}{2\gamma_j u_j^{1/q}}\mathcal{I}\left(\gamma_j>\tfrac{|\theta_j|}{u_j^{1/q}}\right)\\
&\propto 
\exp\{-\lambda_j\gamma_j^q\}\gamma_j^{qa-1-1}\mathcal{I}\left(\gamma_j>\tfrac{|\theta_j|}{u_j^{1/q}}\right) \\
& = \mbox{TGG} \left(\lambda_j^{-1/q} , qa-1,q, \tfrac{|\theta_j|}{u_j^{1/q}}\right)
\end{align*}
where TGG$(a,d,p,x_0)$  denotes the Generalised Gamma distribution left truncated in $x_0$.
\item[Update $\lambda_j$:] For $j=1,\ldots,p$
\begin{align*}
\pi(\lambda_j\mid\text{rest})&
\propto \lambda_j^a\exp\{-\lambda_j\gamma_j^q\}\lambda_j^{b-1}\exp\{-\lambda_j\xi\}\\
&\propto \lambda_j^{a+b-1}\exp\{-\lambda_j(\gamma_j^q+\xi)\}\\
&=\mbox{Gamma}\left(a+b,\gamma_j^q+\xi\right)
\end{align*}
\item[Update $u_j$:] For $j=1,\ldots,p$
\begin{align*}
\pi(u_j\mid\text{rest})&\propto u_j^{1/q}\exp\{-\phi(q)^qu_j\}\tfrac{1}{2\gamma_j u_j^{1/q}}
\mathcal{I}\left(|\theta_j|<\gamma_j u_j^{1/q}\right)\\
&\propto   \exp\{-\phi(q)^qu_j\}
\mathcal{I}\left(u_j>\tfrac{|\theta_j|^q}{\gamma_j^q}\right)  \\
& = \mbox{Truncated Exp} \left(\phi(q)^q,\tfrac{|\theta_j|^q}{\gamma_j^q}\right)
\end{align*}
\item[Update $\bbeta$:]
\[
\pi(\bbeta\mid\text{rest})  = \mbox{N}_p(\hat{\bbeta},\sigma^2(\bX'\bX)^{-1})\prod_{j=1}^p
\mathcal{I}\left(|\theta_j|<\gamma_ju_j^{1/q}\right)
\]
where $\hat{\bbeta}:=(\bX'\bX)^{-1}\bX'\by$ is the ordinary least squares estimate. For generating draws from a truncated multivariate normal distribution, see \cite{rodriguez2004efficient}.
\item[Update $\sigma^2$:]
\begin{align*}
\pi(\sigma^2\mid\text{rest})&\propto 
\left(\tfrac{1}{\sigma^2}\right)^{\alpha+1}\exp\left\{-\tfrac{\theta}{\sigma^2}\right\}
\left(\tfrac{1}{\sigma^2}\right)^{n/2}\exp\left\{-\tfrac{1}{2\sigma^2}  (\by-\bX\bbeta)'   (\by-\bX\bbeta)\right\} \\
& = \mbox{IG}(\alpha+n/2,\theta+ (\by-\bX\bbeta)'   (\by-\bX\bbeta) /2)
\end{align*}
\end{description}

\subsection{MCMC for the case $p \geq n$}

The full conditionals of the random variables in the hierarchical model in Section \ref{sec:MCMC2} are given next.

\begin{description}

\item[Update $\gamma_j$:] 

For $j=1,\ldots,p$
\begin{align*}
\pi(\gamma_j\mid \text{rest}) 
& \propto
\tfrac{1}{\gamma_j}
\exp\left\{ -\tfrac{(\phi(q)|\theta_j|)^q}{\gamma_j^q} - \lambda_j\gamma_j^q \right\}
\gamma_j^{q\tilde{a}-1}\\
 & \propto 
\gamma_j^{q\tilde{a}-2}
\exp\left\{ -\tfrac{(\phi(q)|\theta_j|)^q}{\gamma_j^q} - \lambda_j\gamma_j^q \right\}\\
 & \propto 
\tfrac{1}{q}y_j^{\tilde{a}-1/q - 1}
\exp\left\{ -\left(\tfrac{(\phi(q)|\theta_j|)^q}{y_j} + \lambda_j y_j \right) \right\} 
\quad[\text{Subst } y_j:=\gamma_j^q]\\
& =\textrm{GIG}(2\lambda_j,2(\phi(q)|\theta_j|)^q,\tilde{a}-1/q)
\end{align*}

Here, we considered the representation  $\theta|\kappa \sim \EP(q,\kappa = \frac{1}{1+\gamma^q})$,  $\gamma|\lambda \sim \GG(\lambda^{-1/q},qa,q)$, following Proposition \ref{Thm:fullconditionals}.

\item[Update $\tau_j$:] For $j=1,\ldots,p$

\begin{align*}
\pi(\tau_j\mid \text{rest}) & \propto \tau_j^{-1/2} p_{q/2}(\tau_j)
\tau_j^{1/2}\exp\left\{-\tfrac{\theta_j^2\tau_j\phi(q)^2}{2\gamma_j^2}\right\}\\
& \propto p_{q/2}(\tau_j) 
\exp\left\{-\tfrac{\theta_j^2\tau_j\phi(q)^2}{2\gamma_j^2}\right\}
\end{align*}
This density is an exponentially  tilted $\alpha$-stable distribution 
\citep{devroye2009random} with tilting parameter $\frac{\theta_j^2\phi(q)^2}{2\gamma_j^2}$. 
Random variable generation for exponentially tilted $\alpha$-stable distributions could be done by applying the algorithm in   \cite{devroye2009random} who develops an exact sampling method, which is uniformly fast over all choices of $\alpha$ and the tilting parameter. An alternative exact sampling method, the so-called fast rejection sampling, has been proposed by \cite{hofert2011sampling}. In principle, the sampling method by \cite{hofert2011sampling} works for any exponentially tilted density function over the positive real line. See \cite{favaro2015random} for a discussion of different sampling methods. 

We propose an alternative scheme to generate samples from the full conditional of $\tau_j$, which is potentially more efficient. While the $\tfrac{q}{2}$-stable distribution does not have a closed-form representation (except for $q=1$) when the skewness parameter $1$, it can be written as a scale mixture  (see, \cite[][Theorem 6.6]{devroye2008non} and \cite[][Eq(8)]{choy1997hierarchical}),   
\begin{align*}
p_{q/2}(\tau) 
&= \tfrac{q\tau^{2/(q-2)}}{(2-q)}
\int_0^{1} s(u)\exp\{-s(u)\tau^{q/(q-2)}\}\textrm{d}u 
\end{align*}
with Zolatarev's function \citep{zoloterev1966representation,zolotarev1986one}
\[
s(u) = \left(\tfrac{\sin(\pi uq/2)}{\sin(\pi u)}\right)^{2/(q-2)}
\left(\tfrac{\sin((1-q/2) \pi u)}{\sin( \pi uq/2)}\right)
\] Then

\begin{align*}
\pi(\tau_j\mid \text{rest}) & \propto \tau_j^{-1/2}\tau_j^{2/(q-2)}
\exp\{-s(u_j)\tau_j^{q/(q-2)}\}\tau_j^{1/2}\exp\left\{-\tfrac{\theta_j^2\tau_j\phi(q)}{2\gamma_j}\right\}\\
& =\tau_j^{2/(q-2)}
\exp\{-s(u_j)\tau_j^{q/(q-2)}\}\exp\left\{-\tfrac{\theta_j^2\tau_j\phi(q)}{2\gamma_j}\right\}
\end{align*}

This density is an exponential tilting
\citep{devroye2009random} with tilting parameter $\theta_j^2\phi(q)/(2\gamma_j)$
of a 
log-Gompertz distribution (i.e., the inverse of a Weibull distribution)
\cite[][Section 5.6.1]{kleiber2003statistical} with density given by
\[
f(x) = a\theta^ax^{-a-1}\exp\left\{-(\tfrac{x}{\theta})^{-a}\right\} \quad(x>0)
\]
(the case $a=1$ is a special case of the inverse Gamma distribution)
\item[Update $u_j$:] 
\begin{align*}
\pi(u\mid \text{rest}) & \propto s(u)\exp\{-s(u)\tau^{q/(q-2)}\} \\
\pi(v\mid \text{rest}) & \propto 
\tfrac{e^v}{(1+e^v)^2}s\left(\tfrac{e^v}{1+e^v}\right)\exp\left\{-s\left(\tfrac{e^v}{1+e^v}\right)\tau^{q/(q-2)}\right\} \quad[\text{Subst } v:=\text{logit}(u), \,\tfrac{du}{dv}=\tfrac{u}{1+e^v}]
\end{align*}
We can use the ratio-of-uniforms algorithm (as suggested by \cite{choy1997hierarchical})  to generate random variates $v$ and then transform them back to $u$.

The maximum of $f(x) = x e^{-x\tau^{q/(q-2)}}$ is achieved for $x^*= \tau^{q/(2-q)}$, with $f(x^*)=\tau^{q/(2-q)}e^{-1}$. 
Thus
\begin{align*}
a(r)&:=\sup_{u\in(0,1)} \left[\pi(u\mid \text{rest})\right]^{1/(r+1)} = \tau^{\tfrac{q}{(r+1)(2-q)}}e^{-1/(r+1)}\\
b^-(r)&:=-b^+(r)\\
b^+(r)&:=\sup_{u\in(0,1)}u \left[\pi(u\mid \text{rest})\right]^{r/(r+1)} \leq  [a(r)]^r
\end{align*}

\item[Update $\bbeta$:] 
\begin{align*}
\pi(\bbeta\mid\text{rest})&\propto
\mbox{N}(\by; \bX\bbeta,\sigma^2\bI_n)  \mbox{N}(\bbeta; \mathbf{0},\bS) \\
& = \textrm{N} (\bA^{-1}\bX'\by/\sigma^2,\bA^{-1})
\end{align*}
with $\bA:= \bX'\bX/\sigma^2 +\bS^{-1}$ and $\bS:=\phi(q)^{-2}\diag(\gamma_1^2/\tau_1,\ldots,\gamma_p^2/\tau_p)$ \citep[see][17.2]{west2006bayesian}.

\item[Update $\sigma^2$:]
\begin{align*}
\pi(\sigma^2\mid\text{rest})&\propto 
\left(\tfrac{1}{\sigma^2}\right)^{\alpha+1}\exp\left\{-\tfrac{\theta}{\sigma^2}\right\}
\left(\tfrac{1}{\sigma^2}\right)^{n/2}\exp\left\{-\tfrac{1}{2\sigma^2}  (\by-\bX\bbeta)'   (\by-\bX\bbeta)\right\} \\
& = \mbox{IG}(\alpha+n/2,\theta+ (\by-\bX\bbeta)'   (\by-\bX\bbeta) /2)
\end{align*}
\end{description}

\section{JAGS code}\label{sect:JAGS}

We provide the JAGS code used to fit the Bayesian model described in Section \ref{sec:sim_study} of the paper. The code specifies the model structure, priors, and likelihoods used in the simulations. We consider $n$-dimensional data $\mathbf{y}$ with design matrix $\mathbf{X}$, and for the prior on the regression coefficients $\boldsymbol{\theta}$ we followed the hierarchical model in \eqref{sec:MCMC1}. The hyperparameters are \verb|tau, q, a, b| and \verb|xi|, and \verb|eta| is equal to $\phi(q)$.

\begin{verbatim}
model{
    for(i in 1:n){
      y[i]  ~ dnorm(inprod(X[i,],theta),tau) }
      
    for(i in 1:p){
      theta[i]  ~ dunif( - gamma[i]*u[i]^(1/q), gamma[i]*u[i]^(1/q))
      u[i]      ~ dgamma(1+1/q, eta^q)
      gamma[i]  ~ dggamma(a, pow(lambda[i],1/q), q)
      lambda[i] ~ dgamma(b, xi)}
}
\end{verbatim}

\bibliographystyle{chicago}
\bibliography{prior_references}
\end{document}